\documentclass[12pt]{article}
\usepackage{amscd,amsthm}
\usepackage{latexsym}
\usepackage{amsmath}
\usepackage{amssymb}
\numberwithin{equation}{section}
\def\bbZ{{\Bbb Z}}
\def\bbR{{\Bbb R}}

\def\bbE{{\Bbb E}}
\def\bbC{{\Bbb C}}

\def\Im{{\rm Im}}
\def\const{{\rm const}}

\def\tr{{\rm Tr}}
\def\var{{\rm Var}}
\def\id{{\rm Id}}
\def\parr{{\rm Par}}
\def\cov{{\rm Cov}}

\def\a{{\cal A}}
\def\b{{\cal B}}
\def\r{{\cal R}}
\def\p{{\cal P}}

\allowdisplaybreaks

\begin{document}
\title{Determinantal Random Point Fields}
\author{Alexander Soshnikov\\Caltech\\Department of
Mathematics\\Sloan 253--37\\Pasadena, CA  91125, USA\\sashas@gibbs1.caltech.edu\\
and\\University of California, Davis\footnote{permanent address}\\One Shields
Ave.,\\Davis, CA  95616, USA\\soshniko@math.ucdavis.edu}
\date{}
\maketitle
\begin{abstract}
The paper  contains an exposition  of recent as well as sufficiently old
results  on
determinantal random point fields .
We start with some general theorems  including the proofs  of the
necessary and sufficient condition for the existence of determinantal
random point field with Hermitian kernel and a criterion for the weak convergence of its
distribution. In the second section we  proceed with the examples of the determinantal
random fields from  Quantum Mechanics,
Statistical Mechanics, Random Matrix Theory, Probability Theory,
Representation Theory and Ergodic Theory.
In connection with the Theory of Renewal Processes we 
characterize all Hermitian determinantal random point fields in
$\bbR^1 $ and $\bbZ^1 $
with independent identically distributed spacings.  In the third section
we study the translation
invariant determinantal random point fields and prove the mixing property of
any multiplicity and the absolute continuity of the spectra.
In the last  section we discuss the proofs of the the Central Limit Theorem
for the number of particles in the growing box and the Functional Central Limit
Theorem for the empirical distribution function of spacings. 

\end{abstract}

\section{Definition and General Properties of Determinantal Random Point Fields}

Let $E$ be a one-particle space and $X$ a space of finite or countable configurations
of particles in $E$.  In general $E$ can be a separable Hausdorff space, however
for our purposes it is enough to consider
\begin{equation}
E=\prod^m_{j=1}E_j,\text{ where }E_j\cong \bbR^d\text{ (or }\bbZ^d)
\end{equation}
If it is not mentioned specifically otherwise we always assume below $E=\bbR^d$
with the understanding that all resuls can be easily generalized to (1.1).  We 
assume that each configuration $\xi =(x_i), x_i\in E, i\in\bbZ^1$ (or $\bbZ^1_+$ if
$d>1$), is locally finite, that is for every compact $K\subset E\ \#_K(\xi )=\#
(x_i\in K)$ is finite. The particles in $\xi$ are ordered in some natural way,
e.g., $x_i\leq x_{i+1}$ for $d=1$, and if $d>1$ then either $x_i=x_{i+1}$,
or
\begin{equation}
\begin{split}
\vert x_i\vert =\left (\sum^d_{j=1}(x^{(j)}_i)^2\right )^{\frac{1}{2}}<\vert x_{i+1}
\vert =\left (\sum^d_{j=1}(x^{(j)}_{i+1})^2\right )^{\frac{1}{2}}\\
\end{split}
\end{equation}
where $x_i=(x^{(1)}_i,\dots ,x^{(d)}_i)$,
or
$\vert x_i\vert =\vert x_{i+1}\vert\text{ and there exists }1\leq r\leq d$ such that
$x^{(j)}_i\leq x^{(j)}_{i+1}, 1\leq j\leq r-1$, and $x^{(r)}_i<x^{(r)}_{i+1}$

To define a $\sigma$-algebra of measurable subsets of $X$ we first construct 
the so-called cylinder sets.  Let $B\subset E$ be any bounded Borel set and 
$n\geq 0$.  We call
$C^B_n=\{\xi\in X:\#_B(\xi )=n\}$ a cylinder set.  We define $\b$ as a 
$\sigma$-algebra generated by all cylinder sets (i.e., $\b$ is a minimal 
$\sigma$-algebra that contains all $C^B_n$).
\medskip

\noindent{\bf Definition 1}.
A random point field  is a triplet $(X,\b, P)$  \ where $P$ is  a probability
measure on $(X,\b )$.
\medskip

This definition raises a natural question, namely how one can construct such
probability measures.  The corresponding theory was developed by Lenard in [L1--L3]
where a general case of $E$ locally compact Hausdorff space satisying the second
axiom of countability was studied.  If $E=\bbR^d$ or $\bbZ^d$ one can proceed quite
naively by employing Kolmogorov's fundamental theorem from the theory of stochastic
processes ([K]).  Let $t$ and $s$ be two vectors from $E$ with rational coordinates
$t=(t^{(1)},\dots ,t^{(d)}), s=(s^{(1)},\dots ,s^{(d)})$.  We denote an open
rectangle $\{x=(x^{(1)},\dots x^{(d)})\in E: x^{(j)}=t^{(j)}+\theta_j(s^{(j)}-
t^{(j)}), 0<\theta_j<1,j=1,\dots ,d\}$ by $\sqcap_{t,s}$.  Let us denote the
family of finite unions of open, closed or semi-closed rectangles with rational
$t,s$ by $\r$.  Suppose we are able to construct a joint distribution of non-negative
integer-valued random variables $\eta_D, D\in\r$ (that we later identify with
$\#_D$) such that the following finite-additivity condition holds
\begin{equation}
\eta_D=\sum^n_{i=1}\eta_{D_i}\text{ (a.e.)}
\end{equation}
if $D=\bigsqcup^n_{i=1}D_i,\ \  D,D_i\in\r , i=1,\dots n$.  One immediately can replace
(1.3) then by $\sigma$-additivity property
\begin{equation}
\eta_D=\sum^\infty_{i=1}\eta_{D_i}\text{ (a.e.)},
\end{equation}
$D=\displaystyle\bigsqcup^\infty_{i=1}D_i$, $D,D_i\in\r ,i=1,\dots$, (of course
the fact that $\eta_D$ takes only non-negative integers is essential here!).

It is then easy to see that the joint distribution of random variables $\#_D=
\eta_D$, $D\in\r$ with (1.3) (or (1.4) for that matter) uniquely defines
a probability distribution on $(X,\b )$.

Since in many cases it is convenient to define the distribution of random variables
through their moments the following definition appears natural:
\medskip

\noindent{\bf Definition 2}.  
Locally integrable function $\rho_k:E^k\rightarrow
\bbR^1_+$ is called the $k$-point correlation function of the random point field
$(X,\b ,P)$ if for any disjoint bounded Borel subset $A_1,\dots ,A_m$ of $E$
and $k_i\in\bbZ^1_+$, $i=1,\dots m$, $\displaystyle\sum^m_{i=1}k_i=k$
the following identity holds:
\begin{equation}
\bbE\prod^m_{i=1}\frac{(\#_{A_i})!}{(\#_{A_i}-k_i)!}=\int_{A^{k_1}_1\times\dots\times
A^{k_m}_m}\rho_k(x_1,\dots ,x_k)dx_1\dots dx_k
\end{equation}
where by $\bbE$ we denote the mathematical expectation with respect to $P$.
In particular $\rho_1(x)$ is the density of particles, since
$$\bbE\#_A=\int_A\rho_1(x)dx$$
for any bounded Borel $A\subset E$.  In general $\rho_k(x_1,\dots ,x_k)$ has the
following probabilistic interpretation:  let $[x_1,x_i+dx_i]$, $i=1,\dots ,k$
be infinitesimally small boxes around $x_i$, then $\rho_k(x_1,x_2,\dots ,x_k)
dx_1\cdot\ldots\cdot dx_k$ is the probability to find a particle in each of these
boxes.  The problem of existence and uniqueness of a random point field defined
by its correlation functions was studied in [L1--L3].  Not very surprisingly,
Lenard's papers revealed many similarities to the classical moment problem 
([A], [S2]).  In particular the random point field is uniquely defined by
its correlation functions if the distribution of random variables $\{\#_A\}$
is uniquely determined by its moments.  The sufficient condition for the
uniqueness derived in [L1] reads
\begin{equation}
\sum^\infty_{k=0}\left (\frac{1}{(k+j)!}\int_{A^{k+j}}\rho_{k+j}(x_1,\dots ,x_{k+j})
dx_1,\dots dx_{k+j}\right )^{-\frac{1}{k}}=\infty
\end{equation}
for any bounded Borel $A\subset E$ and any integer $j\geq 0$, however we invite
the reader to check that the divergence of the series with $j=0$, namely
\begin{gather}
\sum^\infty_{k=0}\left (\frac{1}{k!}\int_{A^k}\rho_k(x_1,\dots ,x_k)dx_1,\dots
dx_k\right )^{-\frac{1}{k}}=\infty\tag{1.6'}
\end{gather}
implies (1.6) for any $j\geq 0$.  In [L2], [L3] Lenard obtained the necessary and
sufficient condition for the existence of a random point field with the prescribed
correlation functions.
\medskip

\noindent{\bf Theorem 1}. (Lenard)

{\it Locally integrable functions $\rho_k:E^k\rightarrow\bbR^1$,
$k=1,2,\dots$ are the correlation functions of some random point field if
and only if the Symmetry and Positivity Conditions below are satisfied.}
\begin{enumerate}
\item[a)] {\bf Symmetry Condition}

$\rho_k$ is invariant under the action of the symmetric group $S_k$, i.e.,
\begin{equation}
\rho_k(x_{\sigma (1)},\dots ,x_{\sigma (k)})=\rho_k(x_1,\dots ,x_k)
\end{equation}
for any $\sigma\in S_k$.
\item[b)] {\bf Positivity Condition}

For any finite set of measurable bounded functions $\varphi_k:E^k\rightarrow
\bbR^1$, $k=0,1,\dots ,N$ with compact support, such that 
\begin{equation}
\varphi_0+\sum^N_{k=1}\sum_{i_1\ne\dots\ne i_k}\varphi_k(x_{i_1},\dots ,x_{i_k})\geq 0
\end{equation}
for all $\xi =(x_i)\in X$, the next inequality must be valid:
\begin{equation}
\varphi_0+\sum^N_{k=1}\int_{E^k}\varphi_k(x_1,\dots ,x_k)\rho_k(x_1,
\dots ,x_k) dx_1\dots dx_k\geq 0.
\end{equation}
\end{enumerate}

The necessary part of the theorem is quite easy since both conditions have an obvious
probabilistic interpretation.  In particular the Positivity Condition means that the
mathematical expectations of a certain class of non-negative random variables
must be non-negative.  The sufficient part is more elaborate and relies on an 
analogue of the Riesz Representation Theorem and the Riesz--Krein Extension Theorem
(a close relative of the Hahn--Banach Theorem).  It should be noted that Lenard
established his results in a general setting when $E$ is locally compact
Hausdorff space with the second axiom of countability.  

One can obtain a slightly
weaker (but still hopelessly ineffective!) variant of the Positivity Condition
by approximating $\varphi_k$ from above by step functions.  Let $\p_k$ be the class
of polynomials in $k$ variables that take non-negative values on non-negative integers.
Since the polynomials $\displaystyle\{\prod^k_{i=1}\prod^{m_i-1}_{j=0}(x_i-j),
m_i\geq 0\}$ form a linear basis in the vector space of all polynomials in $k$
variables, we can represent any $q(x_1,\dots ,x_k)\in\p_k$ as
\begin{equation}
q(x_1,\dots x_k)=\sum_{m_1,\dots ,m_k\geq 0}a_{m_1,\dots ,m_k}\cdot\prod^k_{i=1}
\prod^{m_i-1}_{j=0}(x_i-j)
\end{equation}
\medskip

\noindent{\bf Positivity Condition$^*$}: For any $q\in\p_k, k\geq 1$, any bounded
Borel sets $A_1,\dots ,A_k\subset E$, the following condition must be satisfied:
\begin{equation}
a_{0,\dots ,0}+\sum_{m\geq 1}\sum_{m_1+\dots +m_k=m}a_{m_1,\dots ,m_k}
\int_{\prod^k_{i=1}A^{m_i}_i}\rho_m(x_1,\dots ,x_m)dx_1\dots dx_m\geq 0.
\end{equation}
Indeed, the l.h.s. at (1.11) is equal to 
\begin{equation}
\begin{split}
\bbE q(\#_{A_1},\dots ,\#_{A_k})&=\bbE \bigl [a_{0,\dots 0}+\sum_{m\geq 1}
\sum_{m_1+\dots +m_k=m}a_{m_1,\dots ,m_k}\cdot\sum_{i_1\ne\dots\ne i_m}\\
&\chi_{A^{m_1}_1\times\dots\times A^{m_k}_k} (x_{i_1},\dots ,x_{i_m})\bigr ]
\end{split}
\end{equation}
One can notice that in a sense the Positivity Condition$^*$ is similar to
the condition on the moments of the integer-valued nonnegative random 
variable.

In our paper we will study a special class of random point
fields introduced by Macchi in [Ma] (see also [DVJ]).  We start with an
integral operator $K: L^2(\bbR^d)\rightarrow L^2(\bbR^d)$ that we assume
to be non-negative and locally trace class.  The last condition means that
for any compact $B\subset\bbR^d$ the operator $K\cdot\chi_B$ is trace class,
where $\chi_B(x)$ is an indicator of $B$.  Therefore we have
\begin{equation}
K\geq 0,\ \ \tr (\chi_B K\cdot\chi_B)<+\infty
\end{equation}
The kernel of $K$ is defined up to a set of measure zero in $\bbR^d\times
\bbR^d$.  For our purposes it is convenient to choose it in such a way
that for any bounded measurable $B$ and any positive integer $n$
\begin{equation}
\tr ((\chi_B K \chi_B))=\int_{B}K(x,x)
dx
\end{equation}
It appears that one can indeed achieve this.  We start with
\medskip

\noindent{\bf Lemma 1}. ([S3], [AvSS] Remark 3.4)
{\it Let $K$ be trace class on $L^2(\bbR^d)$.  Then its integral kernel
may be chosen so that the function $M(x,y)\equiv K(x,x+y)$ is a continuous
function of $y$ with values in $L^1(\bbR^d)$.  Furthermore if $m(y)=\int
M(x,y)dx$, then $\tr K=m(0)=\int K(x,x)dx$.}
\medskip

\noindent{\bf Proof}.  We give the proof only when $K$ is non-negative.
The general case is quite similar.  Let $\{\lambda_j\}_{j\geq 1}$ is the 
set of non-zero eigenvalues of $K$ and $\{\varphi_j\}_{j\geq 1}$ is the set
of the corresponding eigenfunctions.  The canonical form of $K$ (as a
selfadjoint compact operator) is
\begin{equation}
K=\sum_{j\geq 1}\lambda_j\cdot (\varphi_j,\cdot )\cdot\varphi_j
\end{equation}
Fix $y\in\bbR^d$ and consider $\displaystyle M(x,y)=\sum^\infty_{j=1}
\lambda_j\cdot\varphi_j(x)\cdot\overline{\varphi_j(x+y)}$ as a 
function of $x$.  Since $\Vert\varphi_j(\cdot )\cdot\overline{\varphi_j(
\cdot +y)}\Vert_1=\int_{\bbR^d}\vert\varphi_j(x)\cdot\overline{\varphi_j
(x+y)}\vert dx\leq\Vert\varphi_j\Vert_2\cdot\Vert\varphi_j\Vert_2=1$,
the series defining $M(\cdot ,y)$ converges in $L^1(\bbR^d)$ for any
$y$ and $\displaystyle\Vert M(\cdot ,y)\Vert_1\leq\sum^\infty_{j=1}\lambda_j
=\tr K<+\infty$.  If we now consider $K(x,y)\equiv M(x,y-x)$, it is
well defined for a.e. $(x,y)\in\bbR^d\times\bbR^d$ and gives a kernel
for $K$.  The $L^1$-continuity of $M(\cdot ,y)$ follows from $\displaystyle
\Vert\sum_{j\geq 1}\lambda_j\left (\varphi_j(\cdot )\overline{\varphi_j(\cdot
+y_1)}-\varphi_j(\cdot )\cdot\overline{\varphi_j(\cdot +y_2)}\right )
\Vert_1
\leq\sum^N_{j=1}\lambda_j\cdot\Vert\varphi_j\Vert_2\cdot\Vert\varphi_j
(\cdot +y_1)-\varphi_j(\cdot +y_2)\Vert_2+\sum_{j\geq N}\lambda_j$.
Choosing $N$ sufficiently large one can make $\displaystyle\sum_{j\geq N}
\lambda_j<\tfrac{\epsilon}{2}$.  Choosing $y_1$ sufficiently close to
$y_2$ so that for each $1\leq j\leq N\ \ \Vert\varphi_j(\cdot )-\varphi_j
(\cdot +y_2-y_1)\Vert_2\leq\tfrac{\epsilon}{2\cdot\sum^N_{j=1}\lambda_j}$
we have the first term less than $\tfrac{\epsilon}{2}$ as well.
\hfill$\Box$

With the help of Lemma 1 we derive
\medskip

\noindent{\bf Lemma 2}.  {\it Let $K$ be a non-negative locally trace
class operator on $L^2(\bbR^d)$.  Then its integral kernel can be
chosen in such a way that for any bounded measurable $B\subset\bbR^d$ 
the function
$$M_B(x,y)=
(K\cdot\chi_B)(x,x+y)$$
is a continuous function of $y$ with values in $L^1(B)$.  Furthermore,
$$\tr (\chi_B K \chi_B)=\int_{B}K(x,x)
dx $$}
\medskip

\noindent{\bf Proof}.  Let $K_n=\chi_{[-n,n]^d}\cdot K\cdot\chi_{[-n,n]^d}$.
By Lemma 1 one can choose the kernel $K_n(x,y)$, such that $K_n(\cdot ,
\cdot +y)$ is continuous in $L^1([-n,n]^d)$-norm.  We denote $M_n(x,y)=K_n
(x,x+y)$.  Since $K_{n+1}(x,y)=K_n(x,y)$ for almost all $(x,y)\in
[-n,n]^d\times [-n,n]^d$ we conclude that for almost all $\vert y\vert\leq
n\ M_{n+1}(x,y)=M_n(x,y)$ for a.e. $\vert x\vert\leq n-\vert y\vert$.
The $L^1$-continuity of $M_{n+1}(\cdot ,y), M_n(\cdot ,y)$ allows to 
replace ``almost all $\vert y\vert\leq n$" by ``all $\vert y\vert\leq n$''.
Therefore for any $y$ the values of $M_n(x,y)$ eventually agree for
a.e. $x$.  We denote this value by $M(x,y)$.  The function $M(\cdot ,y)$
inherits local $L^1$-continuity from $\{M_n(\cdot ,y)\}$ and
(1.14) follows. It may be worthwhile to note that for any 
positive integer $k$ and bounded measurable
$ B_1\subset\bbR^d , \ldots, B_k\subset \bbR^d$
the kernel
\begin{equation}
M^{(k)}_B(x,y)= \bigl(\underset{\longleftarrow k \text{times}
\longrightarrow}{(K\cdot\chi_{B_1}) \cdot\ldots\cdot (K\cdot\chi_{B_k})}\bigr)
(x,x+y)
\end{equation}
is also a continuous function of y with the values in $L^1(\bbR^d)$ and
\begin{equation}
 Tr (K\cdot
\chi_{B_1})\cdot\ldots\cdot (K\cdot\chi_{B_k})=\int_{B_1\times\dots\times B_k}
K(x_1,x_2)\cdot\ldots\cdot K(x_k,x_1)dx_1\dots dx_k,
\end{equation}
in particular
\begin{equation}
 Tr (K\cdot
\chi_{B})\cdot\ldots\cdot (K\cdot\chi_{B})=\int_{B^k}
K(x_1,x_2)\cdot\ldots\cdot K(x_k,x_1)dx_1\dots dx_k.
\end{equation}
Indeed, for any two versions of the integral 
kernel  and $ k>1$
the expressions in the last two integrals coincide up to a set of 
measure zero. Since we have already proved that
there exists the kernel of 
\begin{equation}
K\cdot\chi_{B_1} \ldots K\cdot\chi_{B_k}
\end{equation}
satisfying the $L^1$- continuity condition from Lemma 2, the same 
condition is satisfied by any variant of the kernel.

\hfill$\Box$
\medskip

\noindent{\bf Definition 3}.
A random point field in $E$ is called determinantal (or fermion) if
its $n$-point correlation functions are given by
\begin{equation}
\rho_n(x_1,\dots ,x_n)=\det \bigl (K(x_i,x_j)\bigr )_{1\leq i\leq n}
\end{equation}
In  the case $E=\bigsqcup^M_{j=1}E_j, E_j\cong
\bbR^d$ the definition takes the following form:  Let $K$ be a trace class operator on
$L^2\underset{\longleftarrow m\text{ times}\longrightarrow}{(\bbR^d)\oplus\dots
\oplus L^2(\bbR^d)}$. Then $K$ has a matrix valued kernel $(K_{rs}(x,y))_{1\leq
r,s\leq m}, x,y\in\bbR^d$.
\medskip

\noindent{\bf Definition 3'}.  
A random point field in $E$ is called determinantal (or fermion) if its $n$-point
correlation functions are given by
\begin{equation}
\rho_n(x_{11}, x_{12},\dots x_{1i_1},\dots ,x_{m1}, x_{m2},\dots ,
x_{mi_m})=\det
(K_{rs}(x_{ri}, x_{sj}))_{\underset
{1\leq j\leq i_s, s=1,\dots ,m}
{1\leq i\leq i_r, r=1,\dots ,m}},
\end{equation}
where $n=i_1+i_2+\dots +i_m,\ x_{ri}\in E_r,\ 1\leq r\leq m,\ 1\leq 
i\leq i_r$.
\medskip

\noindent{\bf Remark 1}.  If the kernel is Hermitian-symmetric then
the non-negativity of $n$-point correlation functions
implies that the kernel $K(x,y)$ is non-negative definite and therefore indeed
$K$ must be a non-negative operator. It should be noted however that there exist determinantal
random point fields corresponding to non-Hermitian kernels (see the remark after (1.36) and the 
examples in  the sections 2.2 and 2.5).

\medskip

\noindent{\bf Remark 2}.  The condition (1.13) is satisfied for all 
continuous non-negative definite kernels (see [GK], section III.10 or [RS],
vol. III, section XI.4).  In general situation when $K(x,x)$ is locally
integrable, non-negative definiteness of $K(x,y)$ implies that $K_B$ is a
Hilbert--Schmidt operator and one can use a theorem of Gohberg--Krein ([GK],
section III.10, theorem 10.1) that claims that a non-negative
Hilbert-Schmidt operator A is trace class iff
\begin{equation}
\overline{\lim_{h\rightarrow 0}}\frac{1}{(2h)^{2d}}\int\prod^d_{j=1}[2h-\vert x^j
-y^j\vert ]_+A(x,y)dxdy<\infty
\end{equation}
where $t_+=\max (t,0),\ x=(x^1,\dots ,x^d),\ y=(y^1,\dots y^d)$, and TrA is
then given by (1.22).

An interesting generalization of determinantal random point fields, so called 
immanantal random point fields (processes) was introduced by Diaconis and 
Evans  in [DE].

  The classical formula of Fredholm (see [S1], Chapter 3)
claims that a trace class operator with a continuous (in a usual sense) kernel
satisfies
\begin{equation}
\tr \bigl (\wedge^n(A)\bigr )=\frac{1}{n!}\int\det \bigl (A(x_i,x_j)\bigr)_{1
\leq i,j\leq n}dx_1,\dots dx_n
\end{equation}
In general the kernel $K(x,y)$ may not be continuous,
however (1.18) and the Lidskii theorem (see e.g., [RS], volume IV, 
section XIII.17
or [S1], Theorem 3.7) imply 
\begin{equation}
\int_{B^n}K(x_1,x_2)\cdot\ldots\cdot K(x_n,x_1)dx_1\dots dx_n=\sum^\infty_{j=1}
\lambda^n_j(K_B),
\end{equation}
\begin{equation}
\tr \bigl (\wedge^n(K_B)\bigr )=\sum_{j_1<\dots <j_n}\lambda_{j_1}(K_B)\cdot\ldots
\cdot\lambda_{j_n}(K_B)
\end{equation}
Combining (1.24) and (1.25) one arrives at
\begin{equation}
\tr \bigl (\wedge^n(K_B)\bigr )=\frac{1}{n!}\int_{B^n}\det\bigl (K(x_i,x_j)\bigr 
)_{1\leq i,j\leq n}dx_1,\dots dx_n
\end{equation}

It follows then from (1.17) that

\begin{equation}
\begin{split}
&\tr \bigl ((K\cdot\chi_{B_1})\wedge\dots\wedge (K\cdot\chi_{B_n})\bigr )=\\
&\quad \frac{1}{n!}
\int\det \bigl (K(x_i,x_j)\cdot\chi_{B_j}(x_j)\bigr 
)_{1\leq i,j\leq n}dx_1\dots dx_n
\end{split}
\end{equation}
\medskip

\noindent{\bf Definition 4}.  Let the kernel $K$ as in Lemma 2.  We say
that it defines a determinantal random point field $(X,B,P)$ if (1.21) holds.
\medskip

\noindent{\bf Theorem 2}. {\it Let $(X,B,P)$ be a determinantal random point
field with the kernel $K$.  For any finite number of disjoint bounded Borel sets
$B_j\subset E,\ j=1,\dots ,n$, the generating function of the probability distribution
of $\#_{B_j}=\#\{x_i\in B_j\}$ is given by
\begin{equation}
\bbE\prod^n_{j=1}z^{\#_{B_j}}_j=\det ({\rm Id}+\chi_{B}\sum^n_{j=1}(z_j-1)\cdot K\cdot
\chi_{B_j})
\end{equation}
}
\medskip

\noindent{\bf Remark 3}.  (1.28) is the equality of two entire functions.  
The r.h.s. of (1.28) is well defined as a Fredholm determinant of a trace
class operator (see e.g., [RS], volume IV, section XIII.17 or [S1], section
3).

Recall that by definition
\begin{equation}
\bbE\prod^n_{j=1}z^{\#_{B_j}}=\sum^\infty_{k_1,\dots ,k_n=0}P(\#_{B_j}=k_j,
j=1,\dots ,n)\cdot\prod^n_{j=1}z^{k_j}_j
\end{equation}
and
\begin{equation}
\begin{split}
&\det\biggl ({\rm Id}+\chi_{B}\sum^n_{j=1}(z_j-1)\cdot K\cdot \chi_{B_j}
\biggr )=1+\sum^\infty_{m=1}
\sum^n_{j_1,\dots j_m=1}\prod^m_{\ell =1}\\
&\quad (z_{j_\ell}-1)\cdot\tr (\chi_{B}\cdot K\cdot\chi_{B_{j_1}}\wedge\cdots\wedge \chi_{B} K\cdot
\chi_{B_{j_m}})
\end{split}
\end{equation}
\medskip

\noindent{\bf Proof of Theorem 2}.  The Taylor expansion of the generating function
near $(z_1,\dots ,z_n)=(1,\dots ,1)$ is given by
\begin{equation}
\bbE\prod^n_{j=1}z_j^{\#_{B_j}}=1+\sum^\infty_{m=1}\sum_{m_1+\dots +m_n=m}
\bbE\prod^n_{j=1}\frac{(\#_{B_j})!}{(\#_{B_j}-m_j)!\cdot (m_j)!}
\cdot\prod^n_{j=1}(z_j-1)^{m_j}
\end{equation}
The radius of convergence of (1.30) is infinite since
\begin{equation}
\tr (K\cdot\chi_{B_{j_1}}\wedge\dots\wedge K\cdot\chi_{B_{j_m}}\leq\frac{1}{m!}\tr
(K\cdot\chi_B)^m,\text{ where }B=\bigsqcup^n_{j=1}B_j.
\end{equation}
Therefore, it is enough to show that the coefficients in the series (1.30),
(1.31) coincide.  The case $n=1$ follows then from
(1.5), (1.21), (1.26).  Using (1.27) instead of (1.26) we prove the case $n\geq 1$
as well.
\medskip

\noindent{\bf Remark 4}.  Theorem 2 is well known in the Theory of Random Point
Fields (see [DVJ], p. 140, exercise 5.4.9) and in the Random Matrix Theory
(see [TW1]).

As we already mentioned above, if an operator $K$ defines a determinantal random
point field it must be non-negative because of the non-negativity of the correlation
functions.  It follows from Theorem 2, formula (1.28) that $K$ must also be bounded
from above by the identity operator, i.e., $K\leq 1$.  Indeed, suppose $\Vert K
\Vert >1$.  Then there exists a bounded Borel $B\subset E$ such that $\Vert K_B\Vert >1+
\tfrac{\Vert K\Vert -1}{2}>1$.  Let $\lambda_1(K_B)\geq\lambda_2(K_B)\geq\lambda_3
(K_B)\geq\dots$ be the eigenvalues of $K_B$ and choose $0<z_0<1$ so that $1+(z_0-1)
\cdot\lambda_1(K_B)=0$.  Then $\bbE z_0^{\#_B}=\sum^\infty_{k=1}P(\#_B=k)z^k_0=\det
({\rm Id}+(z_0-1)\cdot K_B)=$(by Theorem XIII.106 from [RS])$=\prod_{j\geq 1}(1+
(z_0-1)\cdot\lambda_j(K_B))=0$. Therefore $P(\#_B=k)=0$ for any $k$, a contradiction.
On the other side assume $0\leq K\leq 1$ and let (1.28) define what we hope to be
the distribution of non-negative integer-valued random variables $\{\#_B\}$.
\medskip

\noindent{\bf Lemma 3}.  {\it Let $0\leq K\leq 1$ and $K$ be a locally trace class
operator.  Then (1.28) defines the distribution
of non-negative integer-valued random variables $\{\#_B\}$ with the additional
property that for $B=\bigsqcup^n_{i=1}B_i$
\begin{equation}
\#_B=\sum^n_{i=1}\#_{B_i}\ (a.e.).
\end{equation}
}

$\Box$ 

We need to show three things:  first, that (1.28) defines some 
finite-dimensional
distributions; second, that the finite-dimensional distributions satisfy the additivity
property (1.33); and third, that the finite-dimensional distributions are consistent
and therefore we can apply the Kolmogorov's Fundamental Theorem to
prove the existence of the distribution of $\{\#_B\}$.  Since the Fredholm 
determinant in (1.28) is 1 when $z_i=1,\ i=1,\dots ,n$, the first statement would
follow from the non-negativity of the Taylor coefficients of the Fredholm determinant
at $z_i=0, i=1,\dots ,n$.  Consider $0\leq z_i\leq 1, i=1,\dots ,n$ and assume for
a moment $\Vert K\Vert <1$ (the case $\Vert K\Vert =1$ would be treated later by
a limiting argument).  Let $B=\bigsqcup^n_{i=1}B_i$.  Then $\Vert K_B\Vert
<1$ and $({\rm Id}-K_B)^{-1}$ is a bounded linear operator such that
$({\rm Id}-K_B)^{-1}-{\rm Id}= K_B \cdot ({\rm Id}- K_B)^{-1}$
is trace class.  Applying Theorem XIII.p105 from [RS], vol. IV we obtain
\begin{equation}
\begin{split}
&\det (\id +\chi_B \sum^n_{j=1}(z_j-1)\cdot K\cdot\chi_{B_j})=\det\bigl ((\id-K_B)
\cdot (\id +\sum^n_{j=1}z_j\cdot\\
&(\id - K_B)^{-1}\cdot \chi_{B}\cdot K\cdot\chi_{B_j})\bigr )=\det (\id - K_B)\cdot
\det (\id +\sum^n_{j=1}z_j\cdot\\
&(\id - K_B)^{-1}\cdot \chi_B K\cdot\chi_{B_j})=\det (\id - K_B)\cdot
\sum^\infty_{k=1}\tr (\wedge^k (\sum^n_{j=1}z_j\cdot\\
& (\id - K_B)^{-1}\cdot \chi_B K\cdot\chi_{B_j}))=\sum_{k_1,\dots ,k_n\geq 0}
\frac{(k_1+\dots +k_n)!}{k_1!\dots k_n!}\prod^n_{j=1}\\
&z^{k_j}_j\cdot\det (\id -K_B)\cdot\tr
(\wedge^n_{j=1}(\wedge^{k_j}(\chi_{B_j}(\id
-K_B)^{-1}\cdot \chi_B K\cdot\chi_{B_j})))
\end{split}
\end{equation}
One can see from (1.34) that Taylor coefficients are, up to some positive factors, 
the traces of the exterior products of the non-negative operators, and, therefore,
non-negative.  We conclude that (1.28) defines some finite-dimensional
distributions.  Since $\bbE \prod^n_{i=1}z^{\#_{B_i}}=\det (\id +\chi_B \sum^n_{i=1}
(z-1) K\chi_{B_i})=\det (\id +(z-1)\cdot K_B)=\bbE z^{\#_B}$,
we conclude that $\#_B=\sum^n_{i=1}\#_{B_i}$ (a.e.).  The formula (1.28) defined
the finite dimensional distributions of $\#_{B_i}$ for disjoint compact sets. In 
the case of non-empty self-intersections one represents $B_i$ as $\sqcup C_{k_i}$,
where $\{C_k\}$ are disjoint sets, defines distributions of $\#_{C_k}$
and then uses the additivity property
(1.33) to define the distributions of $\#_{B_i}$.  To prove the consistency of the
finite-dimensional distributions we note that (1.33) allows us to check it only 
for the disjoint $B_1,\dots ,B_{n+1}$.  But then it trivially follows from $\det
(\id +\chi_B \sum^n_{j=1}(z_j-1)\cdot K\cdot\chi_{B_j}+ \chi_B (1-1)\cdot
K\cdot\chi_{B_{n+1}})
=\det (\id +\chi_B \sum^n_{j=1}(z_j-1)\cdot K\cdot\chi_{B_j})$.  The case $\Vert K\Vert <1$
is proven.  Now let $\Vert K\Vert =1$.  Denote by $K^{(\epsilon )}:=K\cdot (1-\epsilon
), \epsilon >0$ and $\#^{(\epsilon )}_B$ the random variables corresponding to the
kernel $K^{(\epsilon )}$.  Since $\Vert K^{(\epsilon )}\Vert <1$ the arguments above
establish the result of Lemma 3 for $K^{(\epsilon )}$.  
It is an easy exercise to see that
$\bbE\prod^n_{i=1}z_i^{\#^{(\epsilon )}_{B_i}}=1+\sum^\infty_{m=1}\sum_{m_1+\dots +
m_n=m}\bbE\prod^n_{j=1}\tfrac{(\#_{B^{(\epsilon )}_j})!}{(\#_{B^{(\epsilon )}_j}
-m_j)!\cdot (m_j)!}\cdot\prod^n_{j=1}(z_1-1)^{m_j}$ uniformly converges 
with all derivatives to $\bbE
\prod^n_{i=1}z^{\#_{B_i}}_i$ on compact sets as $\epsilon\rightarrow 0$.

Lemma 3 is proven.\hfill$\qed$
\medskip

The results above prove
\medskip

\noindent{\bf Theorem 3}.  {\it  Hermitian locally trace class operator $K$ on $L^2(E)$
defines a determinantal random point field if and only if  $0\leq K\leq 1$.  If the
corresponding random point field exists it is unique.}
\medskip

\noindent $\Box$ The necessary and sufficient condition for the existence of the
field has been already established.  The uniqueness result easily follows from
the general criterion (1.6') since
$\frac{1}{k!}\int_{A^k}\rho_k(x_1,\dots ,x_k)dx_1\dots dx_k=\tr (\wedge^k (K_A))
\leq\frac{\tr (K_A)^k}{k!}\leq\frac{1}{k!}$\hfill\qed
\medskip

Consider arbitrary bounded Borel set $B\subset E$.  Then $\tr 
(K_B)=\bbE\#_B
<\infty$ and the number of particles in $B$ is finite with probability 1.  
Let us 
write $X=\bigsqcup_{0\leq k<\infty}C^B_k$, where as before $C^B_k=\{\xi\in X:\#_B
(\xi )=k\}$.  We choose a kernel for $\chi_B\cdot K\cdot\chi_B$ in such a way
(see Lemma 1) that $(\chi_B\cdot K\cdot
\chi_B)(x,x+y)=\sum^\infty_{i=1}\lambda_i
(B)\varphi_i(x)\cdot\overline{\varphi_i(x+y)}$ is a continuous function of $y$
in $L^1(B)$ norm.  Assume for a moment that $ K_B<1$.  Then
\begin{equation}
L_B(x,x+y)=\sum^\infty_{i=1}\frac{\lambda_i(B)}{1-\lambda_i(B)}\cdot
\varphi_i(x)\cdot\overline{\varphi_i(x+y)}
\end{equation}
is also a continuous function of $y$ in $L^1(B)$ norm and is a kernel of
$L_B=(\id -K_B)^{-1} K_B$.  Taking $B_j$
in (1.34) infinitesimally small  one concludes that for each $C^B_k$ the 
distribution of $k$ particles $x_1\leq x_2\leq\dots\leq x_k$
in $B$ has a density with respect to the Lebesgue
measure.  Denoting this density by $p_k(x_1,\dots ,x_k)$ we obtain
\begin{equation}
p_k(x_1,\dots ,x_k)=\det (\id -K_B)\cdot\det\bigl (L_B(x_i,x_j)\bigr 
)_{1\leq i,j\leq k}
\end{equation}
(It should be noted that (1.36) may be nonnegative even for non-Hermitian kernel $K$,
it is easy to see that such $K$ still  has nonnegative minors).
It follows from the definition of $k$-point correlation functions that
\begin{equation}
\rho_k(x_1,\dots ,x_k)=\sum^\infty_{j=1}\frac{1}{j!}\int_{B^j}p_{k+j}(x_1,\dots ,x_k,
x_{k+1},\dots ,x_{k+j})dx_{k+1}\dots dx_{k+j}
\end{equation}
The system of equations can be inversed :
\begin{equation}
p_k(x_1,\dots ,x_k)=\sum^\infty_{j=0}\frac{(-1)^j}{j!}\int_{B^j}\rho_{k+j}
(x_1,\dots ,x_k, x_{k+1},\dots ,x_{k+j})dx_{k+1}\cdots dx_{k+j}
\end{equation}
Functions $p_k(x_1,\dots ,x_k)$ are called Janossy probability densities (see
[DVJ], p.~122) or exclusion probability densities (see
[Ma]).  It is easy to check that
\begin{equation}
\sum^\infty_{j=0}\frac{1}{j!}\int_{B_j}p_j(x_1,\dots ,x_j)dx_1\dots dx_j=1
\end{equation}
The r.h.s. of (1.36)
still makes sense when $\Vert K_B\Vert =\lambda_1(B)=1$ (and
therefore $p_k(x_1,\dots ,x_k)$
are properly defined in this case too).  Indeed, $\det (\id -K_B)=
\prod^\infty_{j=1}(1-\lambda_j(B))$ as a function of $\lambda_1$ has a zero of
order 1 at $\lambda_1=1$.  We claim that $\det (L(x_i,x_j))_{1\leq i,j\leq k}$ has 
a pole at $\lambda_1=1$ also of order 1.  To see this we write $L=\tilde L+\overset
{\approx}{L}$, where $\tilde L_{i,j}=\tfrac{\lambda_1(B)}{1-\lambda_1(B)}\cdot
\varphi_1(x_i)\cdot\overline{\varphi_1(x_j)}$, $\overset{\approx}{L}=\sum_{\ell\geq 2}
\tfrac{\lambda_\ell (B)}{1-\lambda_\ell(B)}\cdot
\varphi_\ell (x_i)\cdot\overline
{\varphi_\ell (x_j)}.$  Then $\det (L(x_i,x_j))_{1\leq i,j\leq k}=\wedge^k
(L(x_i,x_j)_{1\leq i,j\leq k})$,
and we use the fact that rank $(\tilde L)=1$.  If 1 is a multiple eigenvalue
of $K\cdot\Box_B$, say $\lambda_1(B)=\lambda_2(B)=\dots =\lambda_m(B)=1>\lambda_{m+1}
(B)$ one defines $\tilde L_{i,j}=\sum^m_{\ell =1}\tfrac{\lambda_\ell (B)}{1-\lambda_\ell
(B)}\varphi_\ell (x_i)\cdot\overline{\varphi_\ell (x_j)}$ and proceeds in a similar
manner.  
\medskip

\noindent{\bf Remark 5}.
Following Macchi, we call a random point field regular if for any Borel
$B\subset E$ satisfying $\#_B<\infty$ (P-a.e.), the generating function $\bbE z^{\#_B}$
is entire.  It follows from our results (see also Theorem 4 below) that any 
determinantal random point field is regular.
\medskip

\noindent{\bf Remark 6}.  In [Ma] (Theorem 12, p. 113) (see also 
[DVJ], p. 138) Macchi essentially claimed that a necessary and sufficient 
condition on the integral operator
$K$, locally trace class, to define a regular fermion (=determinantal in our
notations) random point field is $0\leq K<1$.  As one can see from Theorem 3 above
this condition is sufficient, but not necessary (as we established in Theorem 3,
the necessary and sufficient condition is $0\leq K\leq 1)$.  For completeness
it should be noted that Macchi studied the case of continuous $K(x,y)$ with Tr$K<\infty$.
\medskip

\noindent{\bf Remark 7}.  Formula (1.36) was established in [Ma], p. 113 (see also
[DVJ], p. 138 and [TW1], p. 820).

We finish \S 1 with a few more results of general nature about determinantal random
point fields.
\medskip

\noindent{\bf Theorem 4} {\it
\begin{enumerate}
\item[a)] The probability of the event that the number of all particles is finite
is either 0 or 1, depending on whether $\tr K$ is finite or infinite.

\item[b)] The number of particles is less or equal to $n$ with probability 1 if and
only if $K$ is a finite rank operator with rank $(K)\leq n$.

\item[c)] The number of particles is $n$ with probability 1 if and only if $K$ is
an orthogonal projector with rank $(K)=n$.

\item[d)] For any determinantal random point field with probability 1 no two
particles coincide.

\item[e)] To obtain results of the theorem for $B\subset E$ one has
to replace $K$ by $K_B$ .
\end{enumerate}
}
\medskip

\noindent{\bf Proof of Theorem 4}.
\begin{enumerate}
\item[a)] One direction is obvious.  Indeed, if $\tr K=\bbE\#_E<+\infty$, then
$\#_E<+\infty$ with probability 1.  Let us now assume $\tr K=+\infty$.  Consider a
monotone absorbing family of compact sets $\{B_j\}^\infty_{j=1}$ (i.e., $B_i\subset
B_{i+1}$ and $\bigcup^\infty_{i=1}B_i=E)$.  Then $\tr K_{B_j}\underset{j\rightarrow
\infty}{\longrightarrow}+\infty$.  Fix arbitrary large $N$.  By the construction
of $\{B_j\}$ we have $P(\#_E\leq N)=\lim_{j\rightarrow\infty}P(\#_{B_j}\leq N)$.  But
$P(\#_{B_j}\leq N)\leq 2^N\cdot\bbE 2^{-\#_{B_j}}=2^N\cdot\det (\id -\tfrac{1}{2}
\cdot K_{B_j})\leq 2^N\cdot e^{-\frac{1}{2}\tr (K_{B_j})}
\underset{j\rightarrow\infty}{\longrightarrow}0$.

\item[b)] If rank $(K)=n$, then writing $K(x,y)=\sum^n_{i=1}\lambda_i\cdot\varphi_i
(x)\cdot\overline{\varphi_i(y)}$ (a.e.), and $\rho_n(x_1,\dots ,x_n)=\det (K(x_i,x_j)
)_{1\leq i,j\leq n}$ we observe that $\rho_m(x_1,\dots x_m)=0$ (a.e.) for any
$m>n$.  Therefore $\bbE\#_E\cdot (\#_E-1)\cdot\ldots\cdot (\#_E-n)=\int\rho_{n+1}
(x_1,\dots x_{n+1})dx_1\dots dx_{n+1}=0$ which implies $\#_E\leq n$ with probability 1.

In the opposite direction, if $\#_E\leq n$ (a.e.) we have 
$\int_{B^{n+1}}\rho_{n+1}(x_1,\dots ,\break x_{n+1}) 
dx_1\dots dx_{n+1}=0$ for any bounded Borel $B\subset E$, therefore $\tr (\wedge^{n+1}
\break (K_B))=0$. Since $K\geq 0$ we obtain rank $(K_B)\leq n$ for arbitrary compact $B$,
which implies rank $(K)\leq n$.

\item[c)] follows from b) and the formula Var$(\#_E)=\tr (K-K^2)=\prod^n_{i=1}\lambda_i
\cdot (1-\lambda_i)$.

\item[d)] Let $B_n=[-n,n]^d$.  It is enough to show that for any $n$ with probability
1 no two particles in $B_n$ coincide.  Let $\epsilon$ be arbitrary small.  Then
$P\{\exists i\ne j:x_i=x_j\in B_n\}\leq P\{\exists i\ne j: \vert x_i-x_j\vert <
\epsilon, x_i\in B_n, x_j\in B_n\}\leq\int_{B_n} (\int_{\vert x-y\vert <\epsilon}
\rho_2(x,y)dx)dy$.  Since $\rho_2(x,y)$ is locally integrable, the last integral
can be made arbitrary small by letting $\epsilon\rightarrow 0$.\hfill$\Box$
\end{enumerate}

The next result gives a criterion for the weak convergence of determinantal
random point fields.
\medskip

\noindent{\bf Theorem 5}.  {\it Let $P$ and $P_n, n=1,2,\dots$ be probability measures
on $(X,B)$ corresponding to the determinantal random point fields defined by
the Hermitian kernels $K$ and $K_n$.  Let $K_n$ converge to $K$ in the weak operator topology and
$\tr (\chi_B K_n \chi_B)\underset{n\rightarrow\infty}{\longrightarrow}\tr (\chi_B K
\chi_B)$ for any bounded Borel $B\subset E$. Then the probability measures $P_n$ converge
to $P$ weakly on the cylinder sets.}
\medskip

\noindent{\bf Proof of Theorem 5}.  It follows from [S1], Theorem 2.20, p. 40 that
the assumptions of the theorem imply 
\begin{equation}
\tr\vert (K_n-K)_B\vert =\Vert (K_n-K)_B\Vert_1\underset{n\rightarrow\infty}{\longrightarrow}0.
\end{equation}
As a consequence of (1.40) we have
\begin{equation}
\tr\cdot (K_n\cdot\chi_{B_1}\cdot\ldots\cdot K_n\cdot\chi_{B_m})\underset{n\rightarrow
\infty}{\longrightarrow}\tr (K\cdot\chi_{B_1}\cdot\ldots\cdot
K\cdot\chi_{B_m})
\end{equation}
for any compact $B_1,\dots ,B_m$.

Thus using (1.26), (1.27) one can see that the joint moments of $\{\#_B\}$
with respect to $P_n$ converge to the joint moments with respect to $P$.  Since the 
moments of $\#_B$ in the case of the determinantal random points define the
distribution of $\#_B$ uniquely one can see (exercise) that $P_n\overset{W}
{\longrightarrow}P$.\hfill$\Box$

The rest of the notes is organized as follows.  Section 2 is devoted to the various
examples of determinantal random point fields arising in Quantum Mechanics, Statistical
Mechanics, Random Matrix Theory, Representation Theory, Probability Theory (Renewal
Process, 2D Random Growth Models).  
In \S 3 we discuss ergodic properties of the translation
invariant determinantal random point fields.  We also point
out a special role played by the sine kernel $K(x,y)=\tfrac{\sin\pi (x-y)}{\pi
(x-y)}$.  In \S 4 we discuss the Central Limit Theorem for the counting measure and the
Functional Central Limit Theorem for the empirical distribution function of spacings.

It is a great pleasure to thank Ya. Sinai for the encouragement to write this 
paper,
B. Simon for the explaination of  the result of Lemma 1 , G. Olshanski for many valuable remarks,
and A. Borodin,
B. Khoruzhenko, R. Killip, and Yu. Kondratiev for useful conversations.

\section{Examples of Determinantal Random Point Fields}

\subsection{Fermion Gas}

Let $H=-\tfrac{d^2}{dx^2}+V(x)$ be a Schr\"odinger operator with discrete spectrum
acting on $L^2(E)$.  Let $\{\varphi_\ell\}^\infty_{\ell =0}$ be an 
orthonormal basis of the eigenfunctions, $H\varphi_\ell =\lambda_\ell\cdot\varphi_\ell,
\ \lambda_0<\lambda_1\leq\lambda_2\leq\dots$.  Consider the $n^{th}$ exterior power
of $H$, $\wedge^n(H):\wedge^n(L^2(E))\rightarrow\wedge^n(L^2(E))$, where $\wedge^n
(L^2(E))=A_nL^2(E^n))$ is the space of square-integrable antisymmetric functions
of $n$ variables and $\wedge^n(H)=\sum^n_{i=1}(-\tfrac{d^2}{dx^2_i}+V(x_i))$.
In Quantum Mechanics $\wedge^n(H)$  describes the Fermi gas with $n$ particles.  
The ground state of the 
Fermi gas is given by
\begin{equation}
\begin{split}
&\psi (x_1,\dots ,x_n)=\\
&\quad\frac{1}{\sqrt{n!}}\sum_{\sigma\in S_n}(-1)^\sigma
\prod^n_{i=1}\varphi_{i-1}(x_{\sigma (i)})=\frac{1}{\sqrt{n!}}\det (\varphi_{i-1}
(x_j))_{1\leq i,j\leq n}
\end{split}
\end{equation}
It could be noted that $\psi (x_1,\dots ,x_n)$ coincides up to a sign
$\epsilon (x_1,\dots ,x_n)$ with the ground state of $\sum^n_{i=1}(-\tfrac{d^2}
{dx^2_i}+V(x_i))$ acting on $S_nL^2(E^n)$ with the boundary conditions 
$\psi\vert_{x_i=x_j}=0$.  According to the postulate of Quantum Mechanics the 
absolute value squared of the ground state defines the probability distribution
of $n$ particles.  We write
\begin{equation}
\begin{split}
&p(x_1,\dots ,x_n)=\vert\psi (x_1,\dots ,x_n)\vert^2=\frac{1}{n!}\det\bigl 
(\varphi_{i-1}(x_j)\bigr )_{1\leq i,j\leq n}\\
&\quad\cdot\det\bigl (\overline{\varphi_{j-1}(x_i)}
\bigr )_{1\leq i,j\leq n}=\frac{1}{n!}\det \bigl (
K_n(x_i,x_j)\bigr )_{1\leq i,j\leq n},
\end{split}
\end{equation}
where $K_n(x,y)=\sum^{n-1}_{i=0}\varphi_{i-1}(x)\overline{\varphi_{i-1}(y)}$ 
is the kernel of
the orthogonal projector onto the subspace spanned by the first $n$ eigenfunctions
of $H$.  We claim that (2.2) defines a determinantal random point field.  Indeed,
the $k$-point correlation functions are given by
\begin{equation}
\begin{split}
&\rho^{(n)}_k(x_1,\dots ,x_n)=\frac{n!}{(n-k)!}\int p_n(x_1,\dots ,x_n)dx_{k+1}
\dots dx_n=\\
&\quad \det \bigl (K_n(x_1,x_j)\bigr )_{1\leq i,j\leq k}
\end{split}
\end{equation}
The last equality in (2.3) follows from the general lemma well known in Random
Matrix Theory.
\medskip

\noindent{\bf Lemma 4 [Me], p. 89)}.  {\it Let $(E,d\mu )$ be a measurable space
and a kernel $K:E^2\rightarrow\bbR^1$ satisfy
\begin{equation}
\int_EK(x,y)\cdot K(y,z)d\mu (y)=K(x,z)
\end{equation}
\begin{equation}
\int_EK(x,x) d\mu (x)={\rm const}
\end{equation}
Then
\begin{equation}
\begin{split}
&\int_E\det\bigl (K(x_i,x_j)\bigr )_{1\leq i,j\leq n}d\mu (x_n)=\\
&\quad ({\rm const}-n+1)
\cdot\det\bigl (K(x_i,x_j)\bigr )_{1\leq i,j\leq n-1}.
\end{split}
\end{equation}
}

We shall consider in more detail two special cases of $H$.  The first case is 
the harmonic oscillator
\begin{enumerate}
\item[a)] $H=-\tfrac{d^2}{dx^2}+x^2, E=\bbR^1$.  Then
\begin{equation}
\varphi_\ell (x)=\frac{(-1)^\ell}{\pi^{\frac{1}{4}}\cdot (2^\ell\cdot\ell !)^{\frac
{1}{2}}}\exp\left (\frac{x^2}{2}\right )\frac{d^\ell}{dx^\ell}\bigl (\exp (-x^2)\bigr )
\end{equation}
are known as Weber-Hermite functions.  To pass to the thermodynamic limit $n\rightarrow
\infty$ we make a proper rescaling
\begin{equation}
x_i=\frac{\pi}{(2n)^{\frac{1}{2}}}y_i,\ i=1,\dots ,n.
\end{equation}
Then the Christoffel-Darboux formula and the Plancherel-Rotach asymptotics of the
Hermite polynomials ([E]) imply that
\begin{equation*}
\begin{split}
&K_n(x_1,x_2)=\sum^{n-1}_{\ell =0}\varphi_\ell (x_1)\varphi_\ell (x_2)=\\
&\quad \left (\frac{n}{2}\right )^{\frac{1}{2}}\left [\frac{\varphi_n(x_1)\cdot
\varphi_{n-1}(x_2)-\varphi_n(x_2)\varphi_{n-1}(x_1)}{x_1-x_2}\right ]
\end{split}
\end{equation*}
has a limit as $n\rightarrow +\infty$
\begin{equation}
K_n(x_1,x_2)\underset{n\rightarrow\infty}{\longrightarrow}K(y_1,y_2)=\frac{\sin\pi
(y_1-y_2)}{\pi (y_1-y_2)}
\end{equation}
The convergence of kernels implies the convergence of $k$-point correlation functions,
which in turn implies the weak convergence of the distribution
$$\left (\frac{\pi}{(2n)^{\frac{1}{2}}}\right )^n p_n\left (\frac{\pi}{(2n)^{\frac
{1}{2}}}y_1,\dots ,\frac{\pi}{(2n)^{\frac
{1}{2}}}y_n\right )dy_1\dots dy_n$$
to the translation-invariant determinantal random point field with the ``sine kernel''
$K(y_1,y_2)=\tfrac{\sin\pi (y_1-y_2)}{\pi (y_1-y_2)}$.

\item[b)] For another example let $E=S^1=\{z=e^{i\theta},0\leq\theta <2\pi\},
H=-\tfrac{d^2}{d\theta^2}$.  Then
\begin{equation}
\begin{split}
&\varphi_\ell (\theta )=\frac{1}{\sqrt{2\pi}}e^{i\ell\theta},\\
&p_n(\theta_1,\dots ,\theta_n)=\frac{1}{n!}\det\left (\sum^{n-1}_{\ell =0}\frac{1}{2\pi}
e^{i\ell (\theta_j-\theta_k)}\right )_{1\leq j,k\leq n}=\\
&\quad\frac{1}{n!}\det \bigl (K_n(\theta_i,\theta_j\bigr )_{1\leq i,j\leq n}
\end{split}
\end{equation}
where
\begin{equation}
K_n(\theta_1,\theta_2)=\frac{1}{2\pi}\frac{\sin\left (\frac{n}{2}\cdot (\theta_2-
\theta_1)\right )}{\sin\left (\frac{\theta_2-\theta_1}{2}\right )}
\end{equation}
After rescaling $\tfrac{n}{2\pi}\theta_i=y_i,\ i=1,\dots ,n$ the rescaled
correlation functions have the same limit as in (2.9), in particular
$$\lim_{n\rightarrow\infty}\frac{2\pi}{n}K_n\left (\frac{2\pi}{n}y_1, \frac{2\pi}{n}
y_2\right )=\frac{\sin\pi (y_2-y_1)}{\pi (y_2-y_1)}.$$
\end{enumerate}

For more information we refer the reader to [D1]-[D3], [L4]-[L5], [Sp].

\subsection{Coulomb Gas at $\beta =2$}

Examples a), b) from \S 2.1 can be reinterpreted as the equilibrium distribution of $n$
unit charges confined to the one-dimensional line (ex. 2.1a)) or 
the unit circle (ex. 2.1b))
repelling each other according to the Coulomb law of two-dimensional electrostatics.
Writing the potential energy as $H(z_1,\dots ,z_n)=-\sum_{1\leq i<j\leq n}\log
\vert z_i-z_j\vert +\sum^n_{i=1}V(z_i)$, where $V$ is an external potential,
we note that the Boltzmann factor 
$\tfrac{1}{Z}\exp (-\beta H(z_1,\dots ,z_n)),
\beta =2$, is exactly $p_n(z_1,\dots z_n)$ in 2.1a) with $V(z)=\tfrac{1}{2}z^2$,
and $p_n(\theta_1,\dots ,\theta_n)$ in \S 2.1b) with $V(z)=0, z_j=e^{i\theta_j},
j=1,\dots ,n$.

The one-component Coulomb gas in two dimensions (a.k.a. a two-dimen-
\break sional one-component
plasma) was studied in a number of papers including [Gin], [Ja1], [Ja2], [AL], [DFGIL],
[FJ1].  This subject is closely related to the theory of non-Hermitian Gaussian
random matrices (to be discussed in \S 2.3d).  The two-component two-dimensional
Coulomb gas (i.e. a system of positively and negatively charged particles) 
was studied
in [Ga], [CJ1]-[CJ3], [AF], [FJ2].  Let us start with a neutral system of
$n$ positive and $n$ negative particles.  After denoting the complex coordinates by
$u_j$ and $v_j, j=1,\dots ,n$, we write the Boltzmann factor at $\beta =2$ as
\begin{equation*}
\begin{split}
&\exp\biggl (2\sum_{1\leq i<j\leq n}(\log\vert u_i-u_j\vert +\log \vert v_i-v_j\vert
-2\log\vert u_i-v_j\vert )\biggr )=\\
&\quad \frac{\prod_{1\leq i<j\leq n}\vert u_i-u_j\vert^2\cdot\vert v_i-v_j\vert^2}
{\prod_{i,j}\vert u_i-v_j\vert^2}=\biggl\vert\det\left (\frac{1}{u_i-v_j}\right 
)_{1\leq i,j\leq n}\biggr\vert^2.
\end{split}
\end{equation*}
Discretizing the model one allows the positive particles to occupy only the sites
of the sublattice $\gamma\cdot \bbZ^2$ and the negative particles to occupy only
the sites of the sublattice $\gamma\cdot (\bbZ^2+(\tfrac{1}{2},\tfrac{1}{2}))$.
The grand canonical ensemble is defined by the partition function (let $\gamma =1$)
\begin{equation*}
\begin{split}
&Z=1+\sum_{u,v}\lambda_+(u)\lambda_-(v)\cdot\frac{1}{\vert u-v\vert^2}+
\left (\frac{1}{2!}\right )^2\\
&\quad\sum_{u_1,u_2,v_1,v_2}\lambda_+(u_1)\lambda_+(u_2)\lambda_+(v_1)
\lambda_+(v_2)\cdot\biggl\vert\det\left (\frac{1}{u_i-v_j}\right )_{1
\leq i,j\leq 2}\biggr\vert^2+\cdots ,
\end{split}
\end{equation*}
where $\lambda_+(u)=e^{-V(u)}, \lambda_-(u)=e^{V(u)}$, are fugacities and $V$ is an
external potential.  One can rewrite the last formula as
\begin{equation*}
\begin{split}
Z=&\det\biggl (\id +\left (\lambda_+\frac{1+\sigma_z}{2}+\lambda_-
\frac{1-\sigma_z}{2}\right )\\
&\cdot\biggl (\frac{\sigma_x+i\sigma_y}{2}\cdot\frac{1}{z-z'}+\frac{\sigma_x-i\sigma_y}
{2}\cdot\frac{1}{\bar z-\overline{z'}}\biggr )\biggr ),
\end{split}
\end{equation*}
where $\sigma_x,\sigma_y,\sigma_z$ are $2\times 2$ Pauli matrices.  In particular
we see that the grand canonical ensemble is a discrete fermion random point field
(the appearance of matrix-valued kernel reflects the fact that $E=\bbZ^2\bigsqcup
(\bbZ^2+(\tfrac{1}{2},\tfrac{1}{2})))$.  Passing to the continuous 
limit $(\gamma =0)$ one
can see that two- and higher order correlation functions have a limit, and the 
limiting kernel $K$ can be expressed in terms of the Green function of a differential
Dirac operator, namely
\begin{equation*}
\begin{split}
K=&\left (m_+\cdot\frac{1+\sigma_z}{2}+m_-\cdot\frac{1-\sigma_z}{2}\right )\cdot\\
&\left (\sigma_x\partial_x+\sigma_y\partial_y+m_+\cdot\frac{1+\sigma_z}{2}+
m_-\frac{1-\sigma_z}{2}\right )^{-1},
\end{split}
\end{equation*}
where $m_+, m_-$ are rescaled fugacities.  In the special case 
$m_+=m_-\equiv$ const
(i.e., $V\equiv 0$), $K=\begin{pmatrix} K++, & K+-\\ K-+, & K--\end{pmatrix}$ can be
expressed in terms of modified Bessel function (for the details see e.g., [CJ3]).

\subsection{Random Matrix Models}

a) Unitary Invariant Ensembles of Hermitian Randon Matrices
\medskip

The probability distribution in \S 2.1a) (formulas (2.2), (2.7)) allows yet another
interpretation.  It is well known in Random Matrix Theory as the distribution of the
eigenvalues in the Gaussian Unitary Ensemble (G.U.E.).
We recall the definition of G.U.E.  Consider the space of $n\times n$ Hermitian 
matrices $\{A=(A_{ij})_{1\leq i,j\leq n}, {\rm Re}(A_{ij})={\rm Re}(A_{ji}), \Im
(A_{ij})=-{\rm Im}(A_{ji})\}$. A G.U.E. random matrix is defined by its probability
distribution
\begin{equation}
P(dA)=\const_n\cdot\exp (-\tr A^2)dA,
\end{equation}
where $dA$ is a flat (Lebesgue) measure, i.e., $dA=\prod_{i<j}d{\rm Re}(A_{ij})\break
d\Im (A_{ij})\prod^n_{k=1}dA_{kk}$.  The definition of G.U.E. is equivalent to the
requirement that $\{{\rm Re}(A_{ij}),\Im (A_{ij}), 1\leq i<j\leq n, A_{kk}, 1\leq k
\leq n\}$ are mutually independent and Re$(A_{ij})\sim N(0,\tfrac{1}{4}), \Im
(A_{ij})\sim N(0,\tfrac{1}{4}), A_{kk}\sim N(0,\tfrac{1}{2})$.  The eigenvalues of
a random Hermitian matrix are real random variables.  For the derivation of their
joint distribution we refer the reader to [De], sections 5.3-5.4 and [Me], chapters
3, 5.  It appears that the density of the joint distribution with respect to the
Lebesgue measure is given exactly by (2.2), (2.7).

We remark that the distribution of a G.U.E. random matrix is invariant under the
unitary transformation $A\rightarrow UAU^{-1}, U\in U(n)$.  A natural generalization
of (2.12) that preserves the unitary invariance is
\begin{equation}
P(dA)=\const_n\cdot\exp (-2\cdot\tr V(A))dA
\end{equation}
where $V(x)$ can be, for example, a polynomial of even degree with a positive leading
coefficients (see [De], section 5).  The derivation of the formula for the joint
distribution of the eigenvalues is very similar to the G.U.E. case.  The density
$p_n(\lambda_1,\dots ,\lambda_n)$ is given by (2.2), where $\{\varphi_\ell (x)\cdot
e^{V(x)}\}^{n-1}_{\ell =0}$ are the first $n$ orthonormal polynomials with respect
to the weight $\exp (-2V(x))$.  $K_n(x,y)$ is then again a kernel of a projector
and therefore satisfies the conditions of Lemma 4. 
\medskip

\noindent b) Random Unitary Matrices
\medskip

Let us consider the group of $n\times n$ unitary matrices $U(n)$.  There exists a
unique translation invariant probability measure on $U(n)$ (see [We]).  It is called
the Haar measure, we will denote it by $\mu_{Haar}$.  The probability density of the
induced distribution of the eigenvalues is given by
$$p_n(\theta_1,\dots ,\theta_n)=(2\pi)^{-n}\cdot \frac{1}{n!}\cdot\prod_{1\leq k<\ell
\leq n}\vert e^{i\theta_k}-e^{i\theta_\ell}\vert^2,$$
which coincides with (2.10)-(2.11) (see [Me], ch. 9-10, [D1]-[D3]).
In the last formula we used the notations
$$\lambda_1=e^{i\theta_1},\dots ,\lambda_n=e^{i\theta_n}.$$
If one starts with the probability measure $\const_n\cdot e^{-\tr V(U)}d\mu_{Haar}
(U)$ on the unitary group instead of the Haar measure, and replaces the monomials
$\tfrac{1}{\sqrt{2\pi}}e^{i\ell\theta}$ by $\psi_\ell (\theta)\cdot e^{-\tfrac{1}{2}
V(\theta )}$, where $\{\psi_\ell\}^{n-1}_{\ell =0}$ are the first $n$ orthonormal
polynomials in $e^{i\theta}$ with respect to the weight $e^{-V(\theta )}d\theta$,
one stills arrives at the formula (2.10) for the $k$-point correlation functions.
\medskip

\noindent c) Random Orthogonal and Symplectic Matrices
\medskip

The distribution of the eigenvalues of a random orthogonal or symplectic matrix
(with respect to the Haar measure) also has a form of a determinantial random point
field with a fixed number of particles.  For the convenience of the reader we draw
below the chart of the kernels appearing in the ensembles of random matrices
from the Classical Compact Groups.
\medskip

\begin{tabular}{|l|l|}
\hline &\\
& $K_n(x,y)$\\
\hline &\\
$U(n)$ & $\frac{1}{2\pi}\cdot\frac{\sin \left (\frac{n}{2}\cdot (x-y)\right )}
{\sin\left (\frac{x-y}{2}\right )}; E=[0,2\pi]$\\
\hline &\\
$SO(2n)$ & $\frac{1}{2\pi}\cdot\left (\frac{\sin\left (\frac{2n-1}{2}\cdot 
(x-y)\right )}{\sin\left (\frac{x-y}{2}\right )}+\frac{\sin\left (
\frac{2n-1}{2}\cdot (x+y)\right )}{\sin\left (\frac{x+y}{2}\right )}
\right ); E=[0,\pi]$\\
\hline &\\
$SO(2n+1)$ & $\frac{1}{2\pi}\cdot
\left (\frac{\sin (n\cdot (x-y))}
{\sin\left (\frac{x-y}{2}\right )}-\frac{\sin (n\cdot (x+y))}
{\sin\left (\frac{x+y}{2}\right )}
\right ); E=[0,\pi]$\\
\hline &\\
$Sp(n)$ & $\frac{1}{2\pi}\left (\frac{\sin (\frac{2n+1}{2}\cdot (x-y))}
{\sin\left (\frac{x-y}{2}\right )}-\frac{\sin (\frac{2n+1}{2}\cdot (x+y))}
{\sin\left (\frac{x+y}{2}\right )}
\right ); E=[0,\pi]$\\
\hline &\\
$0_-(2n+2)$ & the same as for $Sp(n)$\\
\hline
\end{tabular}
\medskip

For additional information we refer the reader  to 
[Jo1], [DS], [KS], [So1], [So2], [So3].
\medskip

\noindent d) Complex Non-Hermitian Gaussian Random Matrices
\medskip

In [Gin] Ginibre considered the ensemble of complex non-Hermitian random $n\times n$
matrices where all $2n^2$ parameters $\{{\rm Re}A_{ij},\Im A_{ij}, 1\leq i,j\leq n\}$
are independent Gaussian random variables with zero mean and variance $\tfrac{1}{2}$.
The joint probability distribution of the matrix elements is then given by the
formula
\begin{equation}
\begin{split}
&P(dA)=\const_n\cdot\exp (-\tr (A^*\cdot A))dA,\\
&dA=\prod_{1\leq j,k\leq n} d{\rm Re}A_{jk}\cdot d\Im A_{jk}.
\end{split}
\end{equation}
The equivalent definition of (2.14) is that $A=\tilde A+i\cdot\overset{\approx}{A}$,
where $\tilde A$ and $\overset{\approx}{A}$ are two independent G.U.E. matrices.
The eigenvalues $\lambda_1,\dots ,\lambda_n$ are complex random variables.  It was
shown that their distribution is given by the determinantal random point field
in $\bbR^2$ with a fixed number of particles $(\#=n)$ and the correlation functions
\begin{equation}
\rho^{(n)}_k(z_1,\dots ,z_k)=\det \bigl (K_n(z_j,\overline{z_m})\bigr )_{1\leq j,m\leq n}
\end{equation}
where $K_n(z_1,\overline{z_2})=\tfrac{1}{\pi}\exp (-\tfrac{\vert z_1\vert^2}{2}-
\tfrac{\vert z_2\vert^2}{2})\cdot\sum^{n-1}_{\ell =0}\tfrac{z^\ell_1\overline{z_2}^\ell}
{\ell !}$.  We mention in passing that $K_n(z_1,\overline{z_2})$ converges to the
kernel 
\begin{equation}
K(z_1,\overline{z_2})=\tfrac{1}{\pi}\exp (-\tfrac{\vert z_1\vert^2}{2}-
\tfrac{\vert z_2\vert^2}{2}+z_1\cdot\overline{z_2})
\end{equation}
which defines the limiting random point field.  A generalization of (2.14) was studied
in [Gir1], [Gir2], [SCSS], [FKS1], [FKS2].  Let $A=\tilde A+i\cdot v\cdot\overset
{\approx}{A}$, where $\tilde A$ and $\overset {\approx}{A}$ are, as above, two
independent G.U.E. matrices and $v$ is a real parameter (it is enough to consider
$0\leq v\leq 1$).  Let us introduce a new parameter $\tau =\tfrac{1-v^2}{1+v^2}$.
The distribution of the matrix elements is given by
\begin{equation}
P(dA)=\const_n\cdot\exp \left (-\frac{1}{1-\tau^2}\tr (A^*A-\tau{\rm Re}(A^2))
\right )dA
\end{equation}
It induces the distribution of the eigenvalues
\begin{equation}
\begin{split}
&p_n(z_1,\dots ,z_n)\prod^n_{j=1}dz_jd\overline{z_j}=\const_n\cdot\exp\biggl [
-\frac{1}{1-\tau^2}\cdot\sum^n_{j=1}\\
&\left (\vert z_j\vert^2-\frac{\tau}{2}(z^2_j+\overline{z_j}^2)\right )\biggr ]
\cdot \prod_{j<k}\vert z_j-z_k\vert^2\cdot\prod^n_{j=1}dz_jd\overline{z_j}
\end{split}
\end{equation}
It should be noted that the expression (2.18) also appeared in the papers by DiFranceso
{\it et al.} ([DFGIL]) and Forrester-Jancovici ([FJ1]) as the Boltzmann factor
of the two-dimensional one-component plasma.  For the calculation of the correlation
functions we refer the reader to [DFGIL],[FJ1],[FKS1],[FKS2].  The crucial role there is
played by the orthonormal polynomials in the complex plane with the weight
\begin{equation}
w^2(z)=\exp\biggl [-\frac{1}{1-\tau^2}\left (\vert z\vert^2-\frac{\tau}{2}(z^2+
\bar z^2)\right )\biggr ]
\end{equation}
Such orthonormal polynomials can be expressed in terms of the Hermite polynomials,
\begin{equation}
\psi_\ell (z)=\frac{\tau^{\frac{\ell}{2}}}{\pi^{\frac{1}{2}}\cdot 
(\ell !)^{\frac{1}{2}}\cdot (1-\tau^2)^{\frac{1}{4}}}
H_\ell\left (\frac{z}{\sqrt\tau}\right ), \ell =0,1,\dots
\end{equation}
where $\sum^\infty_{n=0}H_n(z)\cdot\tfrac{t^n}{n!}=\exp (zt-\tfrac{t^2}{2})$.
We remark that if $\tau =0$ (Ginibre case) $\psi_\ell (z)=\tfrac{1}{\pi^{\frac{1}{2}}
\cdot (\ell !)^{\frac{1}{2}}}\cdot z^\ell$.  Formula for the correlation functions
($\tau$ arbitrary) generalzies (2.15):
\begin{equation}
\begin{split}
\rho^{(n)}_k=\det\left (K_n(z_i,\overline{z_j})\right )_{1\leq i,j\leq k},\\
K_n(z_1,\overline{z_2})=w(z_1)w(\overline{z_2})\cdot\sum^{n-1}_{\ell =0}\psi_\ell
(z_1)\psi_\ell (\overline{z_2})
\end{split}
\end{equation}
In the limit $n\rightarrow\infty\ K_n(z,\overline{z_2})$ converges to 
\begin{equation}
\begin{split}
&K(z_1,\overline{z_2})=\lim_{n\rightarrow\infty}K_n(z_1,\overline{z_2})=\\
&\quad\frac{1}{\pi (1-\tau^2)}\exp\biggl (-\frac{1}{1-\tau^2}\left (\frac{\vert
z_1\vert^2}{2}+\frac{\vert z_2\vert^2}{2}-z_1\overline{z_2}
\right )\biggr )
\end{split}
\end{equation}
We reamrk that the last formula differs from (2.16) only by the trivial rescaling
$z\rightarrow z\cdot\sqrt{1-\tau^2}$.  A special regime, called the regime of
weak non-Hermiticity, was discovered for the model (2.17) by Fyodorov, Khoruzhenko and
Sommers in [FKS1], [FKS2].  Let
\begin{eqnarray*}
{\rm Re}(z_1)&=&n^{\frac{1}{2}}\cdot x+n^{-\frac{1}{2}}x_1,\\
{\rm Re}(z_2)&=&n^{\frac{1}{2}}\cdot x+n^{-\frac{1}{2}}x_2,\\
\Im (z_1)&=&n^{-\frac{1}{2}}\cdot y_1,\\
\Im (z_2)&=&n^{-\frac{1}{2}}\cdot y_2,
\end{eqnarray*}
Assume the parameters $x, x_1, x_2, y_1, y_2$ fixed and take the limit $n\rightarrow
\infty$ in such a way that $\lim_{n\rightarrow\infty}n\cdot (1-\tau )=\tfrac{\alpha^2}
{2}$.  Then
\begin{equation}
\begin{split}
&\lim_{n\rightarrow\infty}\frac{1}{n}K_n(z_1,z_2)=\frac{1}{\pi\alpha}\cdot\\
&\exp\left [-\frac{y^2_1+y^2_2}{\alpha^2}+i\cdot x\cdot\frac{(y_1-y_2)}{2}\right ]
\cdot g_\alpha\left (\frac{y_1+y_2}{2}-i\cdot\frac{(x_1-x_2)}{2}\right ),
\end{split}
\end{equation}
where
\begin{equation}
g_\alpha (y)=\int^{\sqrt{1-\frac{x^2}{4}}}_{-\sqrt{1-\frac{x^2}{4}}}\frac{du}{\sqrt
{2\pi}}\exp\left [-\frac{\alpha^2u^2}{2}-2uy\right ].
\end{equation}
(if $x>2$ the limit in (2.23) is equal to zero).  The formulas (2.23), (2.24) define
determinantal random point field in $\bbR^2$, different from (2.16).

\medskip

\noindent e) Positive Hermitian Random Matrices
\medskip

Following Bronk [Br] we define the Laguerre ensemble of positive Hermitian 
$n\times n$ matrices.  Any positive Hermitian matrix $M$ can be written as $M=A^*A$, 
where $A$ is some complex matrix.  The probability distribution of a random matrix $M$
is given by
\begin{equation}
\const_n\cdot\exp (-\tr A^*A)\cdot [\det (A^*A)]^\alpha dA,
\end{equation}
where $dA$ is defined as in (2.14) and $\alpha >-1$ (the values of $\alpha$
of special interest are $\pm\tfrac{1}{2}, 0$).  The induced probability distribution
of the (positive) eigenvalues is given by
\begin{equation}
\const_n\cdot\exp\left (-\sum^n_{i=1}\lambda_i\right )\cdot\prod^n_{i=1}\lambda^\alpha_i\cdot
\prod_{1\leq i<j\leq n}(\lambda_i-\lambda_j)^2d\lambda_1\dots \lambda_n
\end{equation}
Employing the associated Laguerre polynominals
$$L^\alpha_m(x)\equiv\frac{1}{n!}e^xx^{-\alpha}\frac{d^m}{dx^m}(e^{-x}x^{m+\alpha}),
m=0,1,\dots$$
one can rewrite (2.26) as
\begin{equation}
\frac{1}{n!}\det \biggl (K_n(x_i,x_j)\biggr )_{1\leq i,j\leq n}
\end{equation}
where
\begin{equation}
K_n(x,y)=\sum^{n-1}_{\ell =0}\varphi^{(\alpha )}_\ell (x)\cdot
\varphi^{(\alpha )}_\ell (y),
\end{equation}
and $\{\varphi^{(\alpha )}_\ell (x)=(\Gamma (\alpha +1)\cdot 
\binom{n+\alpha}{n})^{-\tfrac{1}{2}}L^\alpha_\ell (x)\}^\infty_{\ell =0}$
is the orthonormal basis with respect to the weight $e^{-x}\cdot x^\alpha$ on
the positive semiaxis.  Once again Lemma 4 allows us to calculate explicitly
$k$-point correlation functions and show that they are given by determinants
of $k\times k$ matrices with the kernel (2.28).
\medskip

\noindent f) Hermitian Matrices Coupled in a Chain
\medskip

Let $A_1,\dots ,A_p$ be complex hermitian random $n\times n$ matrices with
the joint probability density
\begin{equation}
\begin{split}
&\const_n\cdot\exp\biggl [-\tr \bigl (\frac{1}{2}V_1(A_1)+V_2(A_2)+\dots
+V_{p-1}(A_{p-1})+\frac{1}{2}V_p(A_p)\\
&\quad +c_1A_1A_2+c_2A_2A_3+\dots +c_{p-1}A_{p-1}A_p\bigr )\biggr ]
\end{split}
\end{equation}
We denote the eigenvalues of $A_j$ (all real) by $\widetilde{\lambda_j}=
(\lambda_{j1},\dots ,\lambda_{jn}), j=1,\dots ,p$.  The induced
probability densit of the eigenvalues is then equal to
\begin{equation}
\begin{split}
&P_n(\widetilde{\lambda_1},\dots ,\widetilde{\lambda_p})=\const_n\cdot
\biggl [\prod_{1\leq r<s\leq n}(\lambda_{1r}-\lambda_{1s})(\lambda_{pr}-
\lambda_{ps})\bigr ]\cdot\\
&\quad\biggl [\prod^{p-1}_{k=1}\det\bigl [w_k(\lambda_{kr},\lambda_{k+1s})
\bigr ]_{r,s=1,\dots ,n}\biggr ]
\end{split}
\end{equation}
where
\begin{equation}
w_k(x,y)=\exp\left (-\frac{1}{2}V_k(x)-\frac{1}{2}V_{k+1}(y)+c_kxy\right )
\end{equation}
Eynard and Mehta [EM] established that the correlation functions of this
model $\rho_{k_1,\dots ,k_p}(\lambda_{11},\dots ,\lambda_{1k_1};\dots 
;\lambda_{p1},\dots ,\lambda_{pk_p})=\prod^p_{j=1}\tfrac{n!}{(n-k_j)!}\int
p_n(\widetilde{\lambda_1},\dots ,\widetilde{\lambda_p})\cdot
\prod^p_{j=1}\prod^n_{r_j=k_j
+1}d\lambda_{jr_j}$ can be written as a $k\times k$ determinant with
$k=k_1+\dots +k_p$
\begin{equation}
\det\left [K_{ij}(\lambda_{ir},\lambda_{js})\right ]_{r=1,\dots ,k_i;s=1,
\dots k_j;i,j=1,\dots ,p}
\end{equation}
For the exact formulas for the kernels $K_{ij}(x,y)$ we refer the reader
to [EM] (see also [AM]).  We remark that (2.32) defines a determinantal random point
field with one-particle space $E$ being the union of $p$ copies of $\bbR^1$.
\medskip

\noindent g) Universality in Random Matrix Models.  Airy, Bessel and sine
Random Point Fields.
\medskip

We start with a general class of kernels of the form
\begin{equation}
K(x,y)=\frac{\varphi (x)\cdot\psi (y)-\varphi (y)\psi (x)}{x-y}
\end{equation}
where
\begin{equation}
\begin{split}
&m(x)\varphi '(x)=A(x)\varphi (x)+B(x)\psi (x)\\
&m(x)\psi '(x)=-C(x)\varphi (x)-A(x)\psi (x)
\end{split}
\end{equation}
and $m(x), A(x), B(x), C(x)$ are polynomials.  It was shown by Tracy
and Widom ([TW2]) that Fredholm determinants of integral operators with
kernels (2.33)-(2.34) restricted to a finite union of intervals satisfy
certain partial differential equations.  Airy, Bessel and sine kernels are
the special cases of (2.33), (2.34).  To define sine kernel we set
$\varphi (x)\equiv\tfrac{1}{\pi}\sin( \pi x),\ \psi (x)\equiv\varphi '(x)\ 
(m(x)\equiv 1,\ A(x)\equiv 0,\ B(x)\equiv 1,\ C(x)\equiv\pi^2)$.  For the
Airy kernel $\varphi (x)\equiv A_i(x),\ \psi (x)\equiv\varphi '(x)\ (m(x)
\equiv 1,\ A(x)\equiv 0,\ B(x)\equiv 1,\ C(x)\equiv -x)$. For the Bessel
kernel $\varphi (x)\equiv J_\alpha (\sqrt x),\ \psi (x)\equiv x\varphi '(x)\
(m(x)\equiv x,\ A(x)\equiv 0,\ B(x)\equiv 1,\ C(x)\equiv\tfrac{1}{4}
(x-\alpha^2))$.  Here $A_i(x)$ is the Airy function and $J_\alpha (x)$
is the Bessel function of order $\alpha$ (see [E]).  Writing down these
kernels explicitly we have (see [TW2], [TW3], [TW4])
\begin{equation}
K_{\rm sine}(x,y)=\frac{\sin\pi (x-y)}{\pi (x-y)},
\end{equation}

\begin{equation}
\begin{split}
K_{\rm Airy}(x,y)&=\frac{\a_i(x)\cdot\a'_i(y)-\a_i(y)\a'_i(x)}{x-y}\\
&=\int^\infty_0A_i(x+t)\a_i(y+t)dt
\end{split}
\end{equation}

\begin{equation}
\begin{split}
K_{\rm Bessel}(x,y)&=\frac{J_\alpha (\sqrt x)\cdot\sqrt y\cdot J'_\alpha
(\sqrt y)-\sqrt x J'_\alpha (\sqrt{x})\cdot J_\alpha (\sqrt y)}{2\cdot 
(x-y)}\\
&=\frac{\sqrt x\cdot J_{\alpha +1}(\sqrt x)\cdot J_\alpha (\sqrt y)
- J_\alpha (\sqrt x)\sqrt y\cdot J_{\alpha +1}(\sqrt y)}
{2\cdot (x-y)}
\end{split}
\end{equation}
As we already mentioned, sine kernel appears as a scaling limit in the bulk
of the spectrum in G.U.E. ([Me], Chapter 5).  In its turn, Airy kernel
appears as a scaling limit at the edge of the spectrum in G.U.E. and at
the (soft) right edge of the spectrum in the Laguerre ensemble, while
Bessel kernel appears as a scaling limit at the (hard) left edge in the Laguerre ensemble ([F], [TW3], [TW4]).  Universality conjecture in
Random Matrix Theory asserts that such limits should be universal for a
wide class of Hermitian random matrices.  Recently this conjecture was
proven for unitary invariant ensembles (2.13) in the bulk of the spectrum
([PS], [BI], [DKMLVZ], [De]) and for some classes of Wigner matrices in the
bulk of the spectrum [Jo4] and at the edge [So4].

In the next subsection we completely characterize determinantal random
point fields in $\bbR^1 (\bbZ^1)$ with independent identically distributed
spacings.

\subsection{Determinantal random point fields with i.i.d. spacings.  Renewal processes}

We start with some basic facts from the theory of renewal processes (see
e.g., [Fe], [DVJ]).  Let $\{\tau_k\}^\infty_{k=1}$ be independent identically
distributed non-negative random variables and $\tau_0$ another non-negative
random variable independent from $\{\tau_k\}^\infty_{k=1}$ (in general
the distribution of $\tau_0$ will be different).  We define
\begin{equation}
x_k=\sum^k_{j=0}\tau_j.
\end{equation}
This gives us a random configuration $\{x_k\}^\infty_{k=0}$ in $\bbR^1_+$.
In probability theory a random sequence $\{x_k\}^\infty_{k=0}$ is known
as delayed renewal process.  We assume that the distribution of random
variables $\tau_k, k\geq 1$ has density $f(x)$, called interval distribution
density, and a finite mathematical expectation $E\tau_1=\int^\infty_0
xf(x)dx$.  The renewal density is defined than as
\begin{equation}
\begin{split}
u(x)=&\sum^\infty_{k=1}f^{k*}(x)=f(x)+\int^x_0f(x-y)f(y)dy +\\
&\int^x_0\int^{x-y_2}_0f(x-y_1-y_2)f(y_1)f(y_2)dy_1dy_2+\dots
\end{split}
\end{equation}
One can express higher order correlation functions of the renewal process through
its one-point correlation function and the renewal density.  
Indeed (see [DVJ], p. 136) for $t_1\leq t_2\leq\dots
\leq t_k$ and $k>1$ the following formula takes place
\begin{equation}
\rho_k(t_1,\dots ,t_k)=\rho_1(t_1)\cdot u(t_2-t_1)\cdot u(t_3-t_2)\cdot
\ldots\cdot u(t_k-t_{k-1})
\end{equation}
It follows immediately from the above definitions that a random point
field in $\bbR^1_+$ has i.i.d. nearest spacings iff it is a renewal process
(2.38).  To make this process translation-invariant the probability density of $\tau_0$
must be given by
\begin{equation}
\frac{1}{E\tau_1}\int^{+\infty}_x f(t)dt\ \ \text{([DVJ], p. 72, [Fe], section
XI.3)}
\end{equation}
Then one-point correlation function is identically constant, $\rho_1(x)
\equiv\rho >0$, so (2.40) implies that the distribution of the process
is uniquely defined by the renewal density (in particular one can obtain
$\rho$ from $u(x)$ since $\rho =(E\tau_1)^{-1}$ and the Laplace transforms
of $f$ and $u$ are simply related).  Macchi ([Ma]) considered a special
class of translation-invariant renewal processes with the interval
distribution density
\begin{equation}
f(x)=2\rho (1-2\rho\alpha )^{-\frac{1}{2}}e^{-\frac{x}{\alpha}}\cdot
\sinh \biggl ((1-2\rho\alpha)^{\frac{1}{2}}\cdot\left (\frac{x}{\alpha}
\right )\biggr ),
\end{equation}
where
\begin{equation}
2\rho\alpha\leq 1,\ \rho >0,\ \alpha >0,
\end{equation}
and showed that it is a determinantal random point field with the kernel
\begin{equation}
K(x,y)=\rho\cdot\exp (-\vert x-y\vert /\alpha )
\end{equation}
(restrictions (2.43) are exactly $0<K\leq {\rm Id}$).

In the next theorem we classify all delayed renewal processes that are also
determinantal random point fields in $\bbR^1_+$.
\medskip

\noindent{\bf Theorem 6}.  {\it Determinantal random point field in $\bbR^1_+$
with  Hermitian kernel has i.i.d. spacings if and only if  its kernel satisfies the following two conditions
in addition to $0\leq K\leq{\rm Id}$   locally trace class :

a) for almost all $x_1\leq x_2\leq x_3$
\begin{equation}
K(x_1,x_2)\cdot K(x_2,x_3)=K(x_1,x_3)\cdot K(x_2,x_2),
\end{equation}

b) for almost all $x_1\leq x_2$ the function
\begin{equation}
K(x_2,x_2)-\frac{K(x_1,x_2)\cdot K(x_2,x_1)}{K(x_2,x_1)}
\end{equation}
depends only on the difference $x_2-x_1$.  If a determinantal random point
field is both translation-invariant and with i.i.d. spacings, it is given
by (2.42)-(2.44).}
\medskip

\noindent{\bf Remark 8}.  Of course a translation-invariant d.r.p.f. in $\bbR^1_+$
can be extended in a unique way to the translation-invariant d.r.p.f. in  $\bbR^1$.
\medskip

\noindent{\bf Proof of Theorem 6}.  First we prove the ``only if" part of the
theorem.  Suppose that a determinantal random point field with a kernel
$K(x,y)$ is also a delayed renewal process.  From (2.40), $k=2,3$, we obtain
the formula for the renewal density
\begin{equation}
u(y-x)=K(y,y)-\frac{K(x,y)\cdot K(y,x)}{K(x,x)},\ y\geq x,
\end{equation}
and the expression for $\rho_3(x_1,x_2,x_3), x_1\leq x_2\leq x_3$:
\begin{equation}
\begin{split}
&\rho_3(x_1,x_2,x_3)=K(x_1,x_1)\cdot u(x_2-x_1)\cdot u(x_3-x_2)\\
&\quad =K(x_1,x_1)
\cdot \biggl (K(x_2,x_2)-\frac{K(x_1,x_2)\cdot K(x_2,x_1)}{K(x_1,x_1)}
\biggr )\\
&\quad \biggl (K(x_3,x_3)-\frac{K(x_2,x_3)\cdot K(x_3,x_2)}
{K(x_2,x_2)}\biggr )
\end{split}
\end{equation}
Since with probability 1 there are no particles outside $A=\{x: K(x,x)>0\}$, we
can always consider random point field restricted to $A$.

Comparing
\begin{equation}
\rho_3(x_1,x_2,x_3)=\det (K(x_i,x_j))_{1\leq i,j\leq 3}
\end{equation}
with (2.48), we have 
\begin{equation*}
\begin{split}
&K(x_1,x_2)\cdot K(x_2,x_1)\cdot K(x_2,x_3)\cdot K(x_3,x_2)\cdot\frac{1}
{K(x_2,x_2)}=\\
&\quad -K(x_1,x_3)\cdot K(x_3,x_1)\cdot K(x_2,x_2)+K(x_1,x_2)\cdot 
K(x_2,x_3)\\
&\quad \cdot K(x_3,x_1)+K(x_1,x_3)\cdot K(x_3,x_2)\cdot K(x_2,x_1)
\end{split}
\end{equation*}
which is equivalent to
\begin{equation*}
\begin{split}
&\frac{1}{K(x_2,x_2)}\cdot\biggl (K(x_1,x_2)\cdot K(x_2,x_3)-K(x_2,x_2)\cdot
K(x_1,x_3)\biggr )\\
&\quad\cdot\biggl (K(x_3,x_2)\cdot K(x_2,x_1)-K(x_3,x_1)\cdot K(x_2,x_2)
\biggr )=0
\end{split}
\end{equation*}
The third factor in the last equality is a complete conjugate of the second
factor, and we obtain (2.45).  Condition b) of the theorem has been already
established in (2.47).  For the translation-invariant d.r.p.f. the
kernel $K(x,y)$ depends only on the difference, therefore $K(x,y)=\rho\cdot
e^{-\vert x-y\vert /\alpha}\cdot e^{i\beta (x-y)}$, and the unitary
equivalent kernel $e^{-i\beta x}K(x,y)\cdot e^{i\beta y}$ coincides with
(2.44).  Now we turn to the proof of the ``if" part of the theorem.  Once
we are given the kernel satisfying (2.45) and (2.46) the candidate for
the renewal density must obey $u(x_2-x_1)=K(x_2,x_2)-\tfrac{K(x_1,x_2)\cdot
K(x_2,x_1)}{K(x_1,x_1)}$ for almost all $x_1\leq x_2$.  Let $x_1\leq x_2\leq
\dots\leq x_k$.  Our goal is to deduce the algebraic identity
\begin{equation}
\begin{split}
&\det\biggl (K(x_i,x_i)\biggr )_{1\leq i,j\leq k}=K(x_1,x_1)\cdot
\prod^{k-1}_{i=1}\biggl (K(x_{i+1},x_{i+1})\\
&\quad -\frac{K(x_i,x_{i+1})\cdot K(x_{i+1}, x_i)}{K(x_i,x_i)}\biggr )
\end{split}
\end{equation}
from the basic identities between the commuting variables $K(x_i,x_j),\overline
{K(x_i,x_j)}$ satisfying $K(x_i,x_j)\cdot K(x_j,x_\ell )=K(x_i,x_\ell 
)\cdot K(x_j,x_j),
1\leq i\leq j\leq\ell\leq k, K(x_i,x_j)=\overline{K(x_j,x_i)}$. Let us
introduce $a(x)=K(x,x)\cdot K(0,x)^{-1}, b(x)=K(0,x)^{-1}$.  Then for
$i\leq j$
$$K(x_i,x_j)=a(x_i)\cdot b(x_j)^{-1},$$
$$K(x_j,x_i)=\overline{a(x_i)}\cdot\overline{b(x_j)}^{-1}.$$
This allows us to write the determinant as
\begin{equation}
\begin{split}
&\begin{vmatrix}a(x_1)\cdot b(x_1)^{-1},&a(x_1)\cdot b(x_2)^{-1}, &\dots &
,a(x_1)\cdot b(x_n)^{-1}\\
\overline{a(x_1)}\cdot\overline{b(x_2)}^{-1},&a(x_2)\cdot b(x_2)^{-1},&
\dots &,a(x_2)\cdot b(x_n)^{-1}\\
\dots &\dots &\dots &\dots\\
\overline{a(x_1)}\cdot\overline{b(x_n)}^{-1}, &\overline{a(x_2)}\cdot
\overline{b(x_n)}^{-1},&\dots &,a(x_n)\cdot b(x_n)^{-1}
\end{vmatrix}\\
&=a(x_1)\cdot b(x_1)^{-1}\cdot\prod^{n-1}_{i=1}\biggl (b(x_{i+1})^{-1}
\cdot \overline{b(x_{i+1})}^{-1}\cdot\\
&\bigl (a(x_{i+1})\cdot \overline
{b(x_{i+1})}-\overline{a(x_i)}\cdot b(x_i)\bigr )\biggr ),
\end{split}
\end{equation}
which is exactly the r.h.s. of (2.50).  Once we established $\rho_k(x_1,\dots
,x_k)=\rho_1(x_1)\cdot\prod^{n-1}_{i=1}u(x_{i+1}-x_i), x_1\leq x_2\leq\dots
\leq x_k$, the rest of the proof is quite easy.  Let $p_k(x_1,\dots ,x_k)$
be Janossy densities, i.e.
the probability density of the event to have particles at $x_1,\dots 
,x_k$, and no particles in between.  We recall that $p_k(x_1,\dots ,x_k)=
\sum^\infty_{j=1}\tfrac{(-1)^j}{j!}\cdot\int\rho_{k+j}(x_1,\dots ,x_k;
y_{k+1},\dots ,y_{k+j})dy_{k+1}\dots dy_{k+j}$, where the integration
in the $j$th term is over $\underset{\longleftarrow j\text{ times}
\longrightarrow}{(x_1,x_k)\times\dots\times (x_1, x_k)}$.  We claim that
\begin{equation}
p_k(x_1,\dots ,x_k)=\rho_1(x_1)\cdot\prod^{k-1}_{i=1} f(x_{i+1}-x_i)
\end{equation}
where the interval distribution density $f$ and the renewal distribution
density $u$ are related via the convolution equation
\begin{equation}
u=f+u*f
\end{equation}
Theorem 6 is proven.\hfill$\Box$
\medskip

\noindent{\bf Remark 9}.  The analogue of Theorem 6 is valid in the discrete
case and the proof is the same.  One has to replace (2.42) by the solution
of the discrete convolution equation (2.53) with $u(n)=1-\rho\cdot e^{-2\beta
n}, K(n_1,n_2)=\rho\cdot e^{-\beta\vert n_1-n_2\vert},0<\rho\leq 1,\ \beta
>0$, so that 
\begin{equation}
\hat f(t)=\sum^\infty_{n=0}f(n)e^{int}=\frac{(1-\rho )-(e^{-2\beta}-\rho)
\cdot e^{it}}{(2-\rho )-(2e^{-2\beta}-\rho +1)\cdot e^{it}+e^{-2\beta}
e^{2it}}
\end{equation}

\noindent{\bf Remark 10}.  One can consider a generalization of Theorem 6 to the
case when the multiplicative identity (2.45) still holds, but the renewal
density
$$u(x_1,x_2)=K(x_2,x_2)-\frac{K(x_1,x_2)\cdot K(x_2,x_1)}{K(x_1,x_1)}$$
no longer depends only on the difference of $x_1$ and $x_2$.  Such processes
have independent but not necessarily identically distributed spacings
since the spacings distribution $f(x,y)dy$ depends on the position $x$ of
the left particle.  Thus
\begin{equation}
\begin{split}
&u(x_1,x_2)=f(x_1,x_2)+\int^{x_2}_{x_1}f(x_1,y_1)\cdot f(y_1,x_2)dy_1\\
&\quad +\int^{x_2}_{x_1}\int^{y_2}_{x_1}f(x_1,y_1)\cdot f(y_1,y_2)\cdot
f(y_2,x_2)dy_1dy_2+\dots
\end{split}
\end{equation}
where $f(x,y)$ is a one-parameter family of probability densities, such that
$$\text{supp }f(x,\cdot )\subset [x,+\infty],\ f\geq 0, \int f(x,y)dy=1.$$
We recall the inversion formula for the equation (2.55):
\begin{equation}
\begin{split}
&f(x_1,x_2)=u(x_1,x_2)-\int^{x_2}_{x_1}u(x_1,y_1)u(y_1,x_2)dy_1\\
&\quad +\int^{x_2}_{x_1}\int^{y_2}_{x_1}u(x_1,y_1)\cdot u(y_1,y_2)\cdot
u(y_2,x_2)dy_1dy_2-\dots
\end{split}
\end{equation}
Writing $K(x,y)=a(x)b(y)^{-1}, x\leq y$, where $a(x)=\tfrac{K(x,x)}{K(0,x)},
b(y)=\tfrac{1}{K(0,y)}$, and $u(x,y)=\tfrac{1}{\vert b(y)\vert^2}\cdot
(a(y)\overline{b(y)}-\overline{a(x)}b(x))$, one can in principle characterize
through (2.56) the class of corresponding interval densities $u(x,y)$.

\subsection{Plancherel Measure on Partitions and its 
Generali-\\zations--$z$-Measures and Schur Measures}

By a partition of $n=1,2,\dots$ we understand a collection of non-negative
integers $\lambda =(\lambda_1,\dots ,\lambda_m)$ such that $\lambda_1+\dots
+\lambda_m=n$ and $\lambda_1\geq\lambda_2\geq\dots\geq\lambda_m$.  We
shall use a notation $\parr (n)$ for the set of all partitions of $n$.
For the basic facts about partitions we refer the reader to [St], [Fu],
[Mac], [Sa].  In particular recall that each partition $\lambda$ of $n$ (denoted
$\lambda\vdash n$) can be identified with a Young diagram with $\vert
\lambda\vert =n$ boxes.  A partition $\lambda '$ corresponds to the
transposed diagram.  Let $d$ be the number of diagonal boxes in $\lambda$
(i.e., the number of diagonal boxes in the Young diagram corresponding to
$\lambda$).  We define the Frobenius coordinates of $\lambda$ as
$(p_1,\dots ,p_d\vert q_1,\dots ,q_k)$, where $p_j=\lambda_j-j,\ q_j=
\lambda_j'-j,\ j=1,\cdots ,d$.  The importance of partitions in Representation
Theory can be most easily understood from the fact that the elements of $\parr
(n)$ label the irreducible representations of the symmetric group $S_n$
(see e.g., [Sa], [Fu]).  The Plancherel measure $M_n$ on the set $\parr (n)$
of all partitions of $n$ is given by
\begin{equation}
M_n(\lambda )=\frac{(\dim\lambda )^2}{n!},
\end{equation}
where $\dim\lambda$ is the dimension of the corresponding representation of
$S_n$.  The dimension $\dim\lambda$ can be expressed in terms of the
Frobenius coordinates via a determinantal formula
\begin{equation}
\frac{\dim\lambda}{n!}=\det\left [\frac{1}{(p_i+q_j+1)\cdot p_i!\cdot q_i!}
\right ]_{1\leq i,j\leq d}
\end{equation}
where $\vert\lambda\vert =n$ ([Ol], Proposition 2.6, formula (2.7)).  Let
$\parr =\bigsqcup^\infty_{n=0}\parr (n)$.  Consider the measure $M^\theta$
on Par, which in analogy with statistical mechanics can be called the grand
canonical ensemble:
\begin{equation}
\begin{split}
&M^\theta (\lambda )=e^{-\theta}\cdot\frac{\theta^n}{n!}M_n(\lambda ),\\
&\text{ if }\lambda\in\parr (n), n=0,1,2,\dots ,\ 0\leq\theta <\infty .
\end{split}
\end{equation}
$M^\theta$ is also called the poissonization of the measures $M_n$.  It
follows from (2.59) that $\vert\lambda\vert$ is distributed by the
Poisson law with the mean $\theta$, and $\displaystyle M^\theta (\ \bigl\vert\ 
\vert\lambda\vert =n)=M_n$.
In the Frobenius coordinates measures $M^\theta , M_n$ can be viewed as 
random point fields on the lattice
$\bbZ^1$.  Recently Borodin, Okounkov and Olshanski ([BOO]) and,
independently, Johansson [Jo3]) proved that $M^\theta$ is a determinantal
random point field (to be exact in [Jo3] only the restriction of
$M^\theta$ to the first half of the Frobenius
coordinates $(p_1,\dots ,p_{d(\lambda )})$ was studied, and as a result,
only the part of (2.60) corresponding to $xy>0$ was obtained).  To
formulate the results of [BOO], [Jo3] we define the modified Frobenius
coordinates of $\lambda$ by
$${\rm Fr} (\lambda ):=\{p_1+\frac{1}{2},\dots ,p_d+\frac{1}{2},-q_1-
\frac{1}{2},\dots , -q_d-\frac{1}{2}\}.$$
Let $\rho^\theta_k(x_1,\dots ,x_k)$ be the $k$-point correlation function
of $M^\theta$ in the modified Frobenius coordinates, where 
$$\{x_1,\dots ,x_k\}\subset\bbZ^1+\tfrac{1}{2}.$$
Then
$$\rho^\theta_k(x_1,\dots ,x_k)=\det [K(x_i,x_j)]_{1\leq i,j\leq k},$$
where $K$ is a so-called discrete Bessel kernel,
\begin{equation}
K(x,y)=\begin{cases}\sqrt\theta \cdot\frac{J_{\vert x\vert -\frac{1}{2}}
(2\sqrt\theta )\cdot J_{\vert y\vert +\frac{1}{2}}(2\sqrt\theta )-
J_{\vert x\vert +\frac{1}{2}}(2\sqrt\theta )\cdot J_{\vert y\vert
-\frac{1}{2}}(2\sqrt\theta )}
{\vert x\vert -\vert y\vert},\text{ if }x\cdot y>0,\\
\sqrt\theta \cdot \frac{J_{\vert x\vert -\frac{1}{2}}
(2\sqrt\theta )\cdot J_{\vert y\vert -\frac{1}{2}}(2\sqrt\theta )-
J_{\vert x\vert +\frac{1}{2}}(2\sqrt\theta )\cdot J_{\vert y\vert
+\frac{1}{2}}(2\sqrt\theta )}
{ x - y},\text{ if }x\cdot y<0,
\end{cases}
\end{equation}
where $J_x(\cdot )$ is the Bessel function of order $x$.  
We note that the kernel $K(x,y)$ is not Hermitian symmetric, however 
the restriction of this kernel
to the positive and negative semi-axis is Hermitian.
(2.60)
can be seen as a limiting case of a more general theorem obtained by
Borodin and Olshanski for the so-called $z$-measures (see Theorem 3.3 of
[BO1], also [BO2], [BO3], [KOV] and references therein).  Let $z,z'$ be complex
numbers such that either
\begin{equation}
\begin{split}
&z'=\bar{z}\in\bbC\setminus\bbZ\\
&\text{or}\\
&[z]<\min (z,z')\leq\max (z,z')<[z]+1,
\end{split}
\end{equation}
where $z,z'$ real and [~] denotes the integer part.  Let $(x)_j=x\cdot
(x+1)\cdot\ldots\cdot (x+j-1), (x)_0=1$.  Below we introduce a 
2-parametric family of probability measures $M^{(n)}_{z,z'}$ on $\parr (n)$.
These measures take their origin in harmonic analysis on the infinite
symmetric group ([KOV], [Ol]).  By definition
\begin{equation}
\begin{split}
&M^{(n)}_{z,z'}(\lambda )=\frac{(z\cdot z')^{d(\lambda )}}{(z\cdot z')_n}
\cdot\prod^{d(\lambda )}_{i=1}(z+1)_{p_i}\cdot (z'+1)_{p_i}\cdot\\
&\quad (-z+1)_{q_i}\cdot (-z'+1)_{q_i}\cdot\frac{\dim^2\lambda}{\vert
\lambda\vert !}
\end{split}
\end{equation}
The conditions on $z,z'$ stated above are equivalent to the requirement that
$(z)_j\cdot (z')_j$ and $(-z)_j\cdot (-z')_j$ are positive for any $j=1,2,
\dots$  We note that $M^{(n)}_{z,z'}$ converges to the Plancherel measure
$M_n$ if $z,z'\rightarrow\infty$.  The measure $M^{(n)}_{z,z'}$ is
called the $n$-th level $z$-measure.  Consider now the negative binomial
distribution on the non-negative integers
$$(1-\xi )^{z\cdot z'}\cdot\frac{(z\cdot z')_n}{n!}\xi^n,\ n=0,1,\dots$$
where $\xi$ is an additional parameter, $0<\xi <1$.  The corresponding
mixture of the $n$-level $z$-measures defines measure $M_{z,z',\xi}$ on Par.
We remark that $M_{z,z',\xi}$ degenerates into $M^\theta$ if $z,z'\rightarrow
\infty,\xi\rightarrow 0$ in such a way that $z z'\xi\rightarrow
\theta$.  It was shown in [BO1] that in the modified Frobenius coordinates
$M_{z,z',\xi}$ is a determinantal random point field on $\bbZ^1+\tfrac{1}{2}$.
The corresponding kernel can be expressed in terms of the Gauss hypergeometric
function and is called the hypergeometric kernel.  It appears that a number
of familiar kernels can be obtained in terms of the hypergeometric kernel,
in particular Hermite kernel ((2.2), (2.3), (2.7)), Laguerre kernel
((2.2), (2.28)), Meixner kernel ((2.67) below), Charlier kernel.  For
the hierarchy of the degenerations of the hypergeometric kernel we refer the
reader to [BO2] \S 9.  Recently Okounkov [Ok1] showed that the measures
$M_{z,z',\xi}$ are the special case of an infinite parameter family of
probability measures on Par, called the Schur measures, and defined as
\begin{equation}
M(\lambda )=\frac{1}{z}s_\lambda (x)\cdot s_\lambda (y),
\end{equation}
where $s_\lambda$ are the Schur functions (for the definition of the Schur functions 
see  [St] or [Mac]), $x=(x_1,x_2,\dots )$ and $y=(
y_1,y_2,\dots)$ are parameters such that 
\begin{equation}
Z=\sum_{\lambda\in\parr}s_\lambda (x)\cdot s_\lambda (y)=\prod_{i,j}
(1-x_iy_j)^{-1}
\end{equation}
is finite and $\{x_i\}^\infty_{i=1}=\overline{\{y_i\}^\infty_{i=1}}$.
Measures $M_{z,z',\xi}$ formally correspond to
$\sum^\infty_{i=1}x^m_i=\xi^{\frac{m}{2}}\cdot z,\ \sum^\infty_{i=1}
y^m_i=\xi^{\frac{m}{2}}\cdot z', m=1,2\dots$ .To be precise one should consider
the Newton power sums as real parameters and express the Schur functions as polynomials 
in the power sums.  By now the reader probably
would not be very surprised to learn that the Schur measures also can be
considered as determinantal random point fields ([Ok1], Theorems 1,2)!

\subsection{Two-Dimensional Random Growth Model}

As our last example we consider the following two-dimensional random growth
model ([Jo2]).  Let $\{a_{ij}\}_{i,j\geq 1}$ be a family of independent
identically distributed random variables with a geometric law
\begin{equation}
p(a_{ij}=k)=p\cdot q^k,\ k=0,1,2,\dots
\end{equation}
where $0<q<1, p=1-q$.  One may think about (2.65) as the distribution of the
first success time in a series of Bernoulli trials.  We define
\begin{equation}
G(M,N)=\max_\pi\sum_{(i,j)\in\pi}a_{ij},
\end{equation}
where the maximum in (2.66) is considered over all up/right paths $\pi$
from (1,1) to $(M,N)$, in other words over $\pi =\{(i_1,j_1)=(1,1),
(i_2,j_2),(i_3,j_3),\dots ,\break (i_{M+N-1}, j_{M+N-1})=(M,N)\}$, such that
$(i_{k+1},j_{k+1})-(i_k,j_k)\in\{(0,1),(1,0)\}$.  We mention in passing
that distribution of random variables $\{G(M,N)\}$ can be interpreted in
terms of randomly growing Young diagrams and totally asymmetric exclusion
process with discrete time (for the details see [Jo2]).  Without loss
of generality we may assume $M\geq N\geq 1$.  To state explicitly the
connection to the determinantal random point fields we introduce the 
discrete weight $w^q_K(x)=\binom{x+K-1}{x}\cdot q^x, K=M-N+1$, on non-negative
integers $x=0,1,2\dots$  The normalized orthogonal polynomials $\{M_n(x)
\}_{n\geq 0}$ with respect to the weight $w^q_K$ are proportional to the
classical Meixner polynomials ([Ch]).  The kernel
\begin{equation}
K_{M,N}(x,y)=\sum^{N-1}_{j=0}M_j(x) M_j(y)\left (w^q_K(x)
w^q_K(y)\right )^{\frac{1}{2}}
\end{equation}
satisfies the conditions of Lemma 4 with respect to the counting measure
on non-negative integers.  Therefore
\begin{equation}
P_N(x_1,\dots ,x_N)=\frac{1}{N!}\det\biggl (K_{M,N}(x_i,x_j)\biggr )_{1\leq
i,j\leq N}
\end{equation}
defines a discrete determinantal random point field.  It was shown by
Johansson that the distribution of the random variable $G(M,N)$
coincides with the distribution of the right-most particle in (2.68).
After appropriate rescaling in the limit $N\rightarrow\infty, M\rightarrow
\infty, \tfrac{M}{N}\rightarrow\const$, this distribution 
converges to the distribution of the right-most
particle in the Airy random point field (2.36).
Additional information on the subject of the last two subsections can be found in
the recent papers/preprints [AD], [BDJ1], [BDJ2], [BR1], [BR2], [BR3], [Bor],
[ITW], [Ku], [Ok2], [PS1], [PS2], [TW5], [TW6].

\section{Translation Invariant Determinantal Random Point Fields}

As before $(X,\b ,P)$ denotes a random point field with a one-particle space
$E$, hence $X$ is a space of locally finite configurations of particles
in $E$, $\b$ is a Borel $\sigma$-algebra of measurarble subsets of $X$ and
$P$ is a probability measure on $(X, B)$.  Throughout this section we always
assume $E=\bbR^d$ or $\bbZ^d$.  We define a continuous action 
$\{T^t\}_{t\in E}$ of $E$ on $X$ in a natural way:
\begin{equation}
T^t: X\rightarrow X,\ (T^t\xi )_i=(\xi )_i+t.
\end{equation}
\medskip

\noindent{\bf Definition 5}.  {\it Random point field $(X,\b,P)$ is called 
translation invariant if for any $A\in \b$, any $t\in E$
$$P(T^{-t}A)=P(A).$$}

The translation invariance of a random point field implies the invariance of
$k$-point correlation functions:
\begin{equation}
\begin{split}
&\rho_k(x_1+t,\dots ,x_k+t)=\rho_k(x_1,\dots ,x_k),\text{a.e. }\\
& k=1,2,\dots ,
t\in E.
\end{split}
\end{equation}
Conversely, if $\{\rho_k\}$ are invariant under $\{T^t\}$, then there exists
a corresponding random point field which is translation invariant ([L3]).
In particular, if the translation invariant
correlation functions define $P$ uniquely then the
random point field is translation invariant.  In the case of a determinantal
random point field this implies the following criterion:
a determinantal random point field is translation invariant if and only if
the kernel $K$ is translation invariant, i.e., $K(x,y)=K(x-y,0)=: K(x-y)$.
In this section we restrict our attention to the translation invariant
determinantal random point fields.  We are interested in the ergodic
properties of the dynamical system $(X,B,P,\{T^t\})$.  For the convenience
of the reader recall some basic definitions of Ergodic Theory ([CFS]).
\begin{enumerate}
\item[--] A dynamical system is said to be ergodic if the measure $P(A)$
of any invariant set $A$ equals 0 or 1.

\item[--] A dynamical system has the mixing property of multiplicity
$r\geq 1$ if for any functions $f_0, f_1,\dots , f_r\in L^{r+1}(X,\b ,P)$
we have
\begin{equation}
\lim_{t_1,\dots ,t_r\rightarrow\infty}\int_Xf_0(\xi )f_1(T^{t_1}\xi ):
\dots :f_r(T^{t_1+\dots +t_r}\xi )dF=\prod^r_{i=0}\int_Xf_i(\xi )dP
\end{equation}

\item[--] A dynamical system has an absolute continuous spectra if for any
$f\in L^2(X,B,P)$ orthogonal to constants
\begin{equation}
\int_Xf(\xi )\overline{f(T^t\xi )}dP=\int e^{i(t\cdot\lambda )}h_f(\lambda )
d\lambda ;
\end{equation}
where the integration at the r.h.s. of (3.4) is over $\bbR^d$ in the continuous
case and over $[0,2\pi ]^d$ in the discrete case, and $h_f(\lambda )
d\lambda$ is a finite measure absolutely continuous with respect to the 
Lebesgue measure.  One can interpret (3.4) in the following way.  We define
a $d$-parameter group of unitary operators $\{U^t\}_{t\in E}$ on $L^2
(X,B,P)$ as
$$(U^tf)(\xi )=f(T^t\xi ).$$
Usually such family of unitary operators is called adjoint to the
dynamical system.  It is easy to see that $\{U^t\}$ commute.  Since $L^2
(X,\b ,P)$ is separable and $(U^t\psi ,\varphi )$ is a measurable function
of $t$ for any $\psi ,\varphi\in L^2(X,B,P)$ one can apply the von Neumann
theorem ([RS], vol. 1, Theorem VIII.9) to conclude that $U^t$ is strongly
continuous.  In the case $E=\bbR^d$ one has $h_f(\lambda )d\lambda =d(f,
Q_\lambda f)$, where $dQ_\lambda$ is a projection-valued measure, $Q_\lambda =
Q_{(-\infty ,\lambda_1)\times\dots\times (-\infty ,\lambda_d)}=\prod^d_{j=1}
\chi_{(-\infty ,\lambda_j)}(A_j), \{A_j\}_{j=1}^d$ are the generators of the
one-parameter groups $U^{(0,\dots ,t_j,0,\dots 0)}$ and $\chi_{(-\infty ,t)}$
is the indicator of $(-\infty ,t)$ ([RS], vol. I, Theorem VIII.12).  In the
discrete case $E=\bbZ^d\ \ dQ_\lambda$ is a projection-valued measure on a
$d$-dimensional torus,
$$Q_{[1,e^{i\lambda_1}]\times\dots\times [1,e^{i\lambda_d}]}=\prod^d_{j=1}
\chi_{[1,e^{i\lambda_j}]}(U_j), U_j=U^{(0,\dots ,t_j=1,\dots 0)}.$$
\end{enumerate}
\medskip

\noindent{\bf Theorem 7}.  {\it Let $(X,B,P)$ be a translation invariant
determinantal random point field.  Then the dynamical system $(X,B,P,\{T^t\})$
is ergodic, has the mixing property of any multiplicity and its spectra is
absolutely continuous.}
\medskip

\noindent{\bf Remark 11}. Recall that the absolute continuity of the spectra
implies the mixing property of multiplicity 1, which in turn implies
ergodicity ([CFS]).
\medskip

\noindent{\bf Proof of Theorem 7}.  We note that the linear combinations
of 
\begin{equation}
\begin{split}
&f(\xi )=\prod^N_{j=1}S_{g_j}(\xi ),\\
&N\geq 1, S_g(\xi )=\sum_ig(x_i),
g_j\in C^\infty_0(\bbR^d), j=1,\dots N
\end{split}
\end{equation}
are dense in $L^2(X,B,P)$.  Therefore it is enough to establish (3.3), (3.4)
for the functions of such form.  We start with the lemma calculating the
mathematical expectation of (3.5).
\medskip

\noindent{\bf Lemma 5}. {\it

a)
\begin{equation}
\begin{split}
&\bbE_P\prod^N_{j=1}S_{g_j}(\xi )=\sum^N_{m=1}
\sum_{\stackrel{\text{over partitions}}{\bigsqcup^m_{\ell =1}C_\ell =\{1,
\dots ,N\}}}\prod^m_{\ell =1}\biggl [\sum^{\#(C_\ell )}_{k_\ell =1}
\sum_{\stackrel{\text{over partitions}}{\bigsqcup^{k_\ell}_{i=1}
B_{\ell i}=C_\ell}}\biggl \{\sum_{\sigma\in S^{k_\ell}}\\
&\frac{(-1)^\sigma}{k_\ell}\cdot\int\prod^{k_\ell}_{i=1}g_{B_{\ell \sigma(i)}}
(x_{i})\cdot K(x_{i+1}-x_{i})
dx_1\dots dx_{k_\ell}\biggr\}
\biggr ]
\end{split}
\end{equation}
where $g_{B_{\ell i}}(x)=\prod_{j\in B_{\ell i}}g_j(x)$.

b) $\bbE\prod^{N_1+\dots +N_{r+1}}_{j=1}S_{g_j}(\xi )
-\prod^{r+1}_{s=1}\left (\bbE\prod^{N_1+\dots +N_s}_{N_1+\dots +
N_{s-1}+1}S_{g_j}(\xi )\right )=($similar expression to (3.6), with the
only difference that partitions 
\begin{equation}
\bigsqcup^m_{\ell =1}C_\ell =\left\{1,2,\dots ,\sum^{r+1}_{s=1}N_s
\right\}
\end{equation}
satisfy (*)), where

(*) There exists at least one element $C_\ell$ of the partition such that
the intersections of $C_\ell$ with at least two of the following sets $\{
1,\dots N_1\},\dots ,\{N_1+\dots +N_{s-1}+1,\dots ,N_1+\dots +N_s\},
\dots ,\{N_1+\dots +N_r +1,\dots ,N_1+\cdots +N_{r+1}\}$ are
non-empty.}
\medskip

\noindent{\bf Proof of Lemma 5}.  The proof of part a) is rather
straightforward and quite similar to the one given at the beginning of
\S 2 in [So3] (see formulas (2.1)-(2.7) from the reference).  The proof of
part b) follows from a).\hfill$\Box$

To derive the mixing property (3.3) we replace $g_j(\cdot )$ for $N_1+
\dots +N_{s-1}+1\leq j\leq N_1+\dots +N_s, s=1,\dots , r+1$, in (3.7) by
$g_j(\cdot +t_1+\dots +t_{s-1})$.  Fix a partition $\bigsqcup^m_{\ell =1}
C_\ell =\{1,2,\dots ,N_1+\dots +N_{r+1}\}$.  Since $\{g_j\}$ are
bounded functions with compact support, each of $m$ factors at the r.h.s.
of (3.7) is bounded.  We claim that the $\ell^{th}$ factor (corresponding
to $C_\ell$, where $\ell$ is the same index as in (*)) goes to
zero.  To see this we fix an arbitrary partition of $C_\ell , 
\bigsqcup^{k_\ell}_{i=1}B_{\ell i}=C_\ell$.  By assumption, $C_\ell$ contains
indices $1\leq u<v\leq N_1+\dots +N_{r+1}$, such that $u$ and $v$
belong to different subsets $\{1,\dots ,N_1\},\dots \{N_1+\dots +N_{s-1}+1,
\dots ,N_1+\dots +N_s\},\dots ,\{N_1+\dots +N_r +1,\dots ,N_1+\dots
+N_{r+1}\}$.  We claim that
\begin{equation}
\int\prod^{k_\ell}_{i=1}g_{B_{\ell \sigma(i)}}(x_{i})\cdot K(x_{i+1}-
x_{i})dx_1\dots dx_{k_\ell}
\end{equation}
goes to zero as min$\{t_s,1\leq s\leq r\}\rightarrow\infty$.  Indeed if
min$\{t_s,1\leq s\leq r\}$ is sufficiently large, the indices $u,v$ belong
to different $B_{\ell i}$'s or the corresponding $g_{B_{\ell i}}$ is zero
(the supports of the factors in $g_{B_{\ell i}}$ will not intersect).  Once $u$ and
$v$ belong to different $B_{\ell i}$'s the argument in $K(x_{i+1}-
x_{i})$ for some $i$ is greater than min$\{t_s,1\leq s\leq r\}$.
Since the Fourier transform of $K(x), \hat K(t)=\int e^{ixt}K(x)dx$
is a non-negative integrable function (bounded from above by 1), applying the
Riemann-Lebesgue lemma we obtain that $K(x_{i+1}-x_{i})$ goes
to zero.  The other terms in (3.8) are bounded and the integration is over a 
bounded set, therefore (3.8) goes to zero and the proof of the mixing property
follows.

To establish the absolute continuity of the spectrum we apply (3.7) when
$r=2, N_1=N_2=N, g_{N+j}(x)=\overline{g_j(x+t)}, j=1,\dots ,N, f(\xi )=
\prod^N_{j=1}S_{g_j}(\xi ),\break \overline{f(T^t\xi )}=\prod^N_{j=1}S_{\overline{g_j}}
(T^t\xi )=\prod^{2N}_{j=N+1} S_{\overline{g_j}}(\xi )$.  We have
\begin{equation}
\begin{split}
&\bbE (f(\xi )-\bbE f)\cdot (\overline{f(T^t\xi )}-\overline{\bbE f})=\sum^{2N}_{m=1}
\sideset{}{^*}\sum_{\stackrel{\text{over partitions}}{\bigsqcup^m_{\ell =1}C_\ell =\{1,\dots , 2N\}}}
\prod^m_{\ell =1}
\biggl [\sum^{\# (C_\ell )}_{k_\ell =1}\\
&\sum_{\stackrel{\text{over partitions}}
{\bigsqcup^{k_\ell}_{i=1}B_{\ell i}=C_\ell}}\bigl \{\sum_{\sigma\in S^{k_\ell}}
\frac{(-1)^\sigma}{k_\ell}\cdot\int\prod^{k_\ell}_{i=1}g_{B_{\ell\sigma (i)}}
(x_i)\cdot K(x_{i+1}-x_i)dx_1\dots dx_{k_\ell}\bigr \}\biggr ],
\end{split}
\end{equation}
where we assume that $x_{k_\ell +1}=x_1$ in the integral, and the sume in $\sum^*$
is over partitions $\{C_1,\dots ,C_m\}$ such that for at least one element $C_\ell$
of the partition both $C_\ell\cap\{1,2,\dots ,N\}$ and $C_\ell\cap\{N+1,\dots ,2N\}$
are non-empty (we denoted above this property by (*)).  The terms in the product
$\prod^m_{\ell =1}$ corresponding to those $\ell$ that do not satisfy (*) are
constants as functions of $t$.  Fix now $\ell$ satisfying (*).  We claim that
\begin{equation}
\begin{split}
&\int\prod^{k_\ell}_{i=1}g_{B_{\ell\sigma (i)}}(x_i)\cdot K(x_{i+1}-x_i)dx_1\dots
dx_{k_\ell}=\\
&\left (\frac{1}{2\pi}\right )^{k_\ell}\prod^{k_\ell}_{i=1}
\widehat g_{B_{\ell\sigma (i)}}
(y_{i+1}-y_i)\cdot \hat K(y_{i+1})dy_1\dots dy_{k_\ell}
\end{split}
\end{equation}
can be written as $\int e^{i(t\cdot\lambda )}h(\lambda )d\lambda$, where $h(\lambda )$
is an integrable function.  The check is rather straightforward and we leave
the details to the reader.  We infer that (3.9) is a linear combination of the
products of the Fourier transforms of integrable functions.  Since the product of
the Fourier transforms is the Fourier transform of the convolution the proof of
the absolute continuity of the spectrum follows.  Theorem 7 is proven.\hfill$\Box$
\medskip

One can without difficulty calculate the spectral density of the centralized
linear statistics
$$S_g(\xi )-\bbE S_g=\sum_ig(x_i)-\bbE\sum_ig(x_i).$$
Namely
\begin{equation}
\begin{split}
&\bbE (S_g-\bbE S_g)(\overline{S_g(T^t\cdot )}-\bbE\overline{S_g})=\int e^{i(t
\lambda )}\cdot (K(0)-\widehat{\vert K\vert^2}(\lambda ))\\
&\cdot\frac{1}{2\pi} \vert\hat g(\lambda )\vert^2d\lambda ,\text{ and }
h_{S_g}(\lambda )=(K(0)-\widehat{\vert K\vert^2}(\lambda ))\cdot\frac{1}{2\pi}
\vert\hat g(\lambda )\vert^2
\end{split}
\end{equation}
We conclude that
\begin{equation}
\mu (d\lambda )=(K(0)-\widehat{\vert K\vert^2}(\lambda ))d\lambda
\end{equation}
is the spectral measure of the restriction of $\{U^t\}$ to the subspace of the
centralized linear statistics.  Since $0\leq\hat K(\lambda )\leq 1, K(0)=\tfrac
{1}{2\pi}\int\hat K(\lambda )d\lambda$, we see that
$$0\leq\frac{d\mu}{d\lambda}=K(0)-\widehat{\vert K\vert^2}(\lambda )=K(0)-\frac{1}
{2\pi}\int\hat K(y)\hat K(y-\lambda )dy\leq K(0).$$
We note that $\tfrac{d\mu}{d\lambda}>0$ for $\lambda\neq 0$, and $\tfrac
{d\mu}{d\lambda}(0)=0$ if and only if $\hat K(\lambda )$ is an indicator.
In particular the spectral measure $\mu$ is equivalent to the Lebesgue
measure.

Before we formulate the next lemma recall that by $\#_{[-L,L]^d}(\xi )$
we denote the number of particles in $[-L, L]^d$.
\medskip

\noindent{\bf Lemma 6}.
\begin{equation}
{\rm Var} (\#_{[-L,L]^d})={\rm Vol}([-L,L]^d)\cdot \left (\frac{d\mu}
{d\lambda}(0)+\bar o(1)\right )\text{ as }L\rightarrow\infty .
\end{equation}
\medskip

\noindent{\bf Proof of Lemma 6}.  The probabilitists are well familiar with
the analogue of this result in the Theory of Random Processes:  let $\{\eta_n
\}$ be $L^2$-stationary random sequence and $h(\lambda )$ its spectral
density,
$\bbE\eta_n\overline{\eta_m}=b(n-m)=\tfrac{1}{2\pi}\int^{2\pi}_0e^{i\lambda
\cdot (n-m)}h(\lambda )d\lambda$, then Var$(\eta_n+\dots +\eta_n)=(h(0)+
\bar o(1))\cdot n$ ([IL], section XVIII.2).  To prove the lemma we write
\begin{equation*}
\begin{split}
&{\rm Var} (\#_{[-L,L]^d})=\int_{[-L,L]^d}\int_{[-L,L]^d}\rho_2(x,y)-
\rho_1(x)\rho_1(y)dxdy\\
&+\int_{[-L,L]^d}\rho_1(x)dx=-\int_{[-L,L]^d}\int_{[-L,L]^d}
\vert K\vert^2 (x-y)dxdy\\
&+K(0){\rm Vol}([-L,L]^d)=(K(0)-\int_{\bbR^d}\vert K\vert^2(x)dx+\bar o
(1))\cdot{\rm Vol}([-L,L]^d)=\\
&(K(0)-\widehat{\vert K\vert^2}(0)+\bar o(1))\cdot{\rm Vol}([-L,L]^d).
\end{split}
\end{equation*}
\hfill$\Box$

The subleading terms in (3.13) also depend on the behavior of $\tfrac{d\mu}
{d\lambda}$ near the origin.  For example, let $\hat K(\lambda )$ be an 
indicator, $\hat K(\lambda )=\chi_B(\lambda ), B\subset\bbR^d$.  As we have
seen above this is equivalent to $\tfrac{d\mu}{d\lambda}(0)=0$.  For
simplicity we will assume $d=1$.  If $B$ is a union of $m$ disjoint
intervals
\begin{equation}
\begin{split}
&\frac{d\mu}{d\lambda}(\lambda )=K(0)-\frac{1}{2\pi}\int\hat K(y)\cdot\hat K
(y-\lambda )dy=\\
&\frac{1}{2\pi}\cdot\left [{\rm length}(B)-{\rm length}(B\cap (B+\lambda 
))\right ]=\\
&\frac{m}{2\pi}\cdot\vert\lambda\vert\cdot (1+\bar o(1)),\ \lambda
\rightarrow 0,
\end{split}
\end{equation}
and after more careful evaluation of the asymptotics of $\int^L_{-L}
\int^L_{-L}\vert K\vert^2(x-y)dxdy=\tfrac{1}{2\pi}\cdot\int^\infty_{-\infty}
\widehat{\vert K\vert^2}(\lambda )\cdot (\tfrac{2\sin (L\cdot\lambda)}
{\lambda})^2d\lambda$ we arrive at
\begin{equation}
{\rm Var}(\#_{[-L,L]^d})=\frac{m}{\pi^2}\log L\cdot (1+\bar o(1))
\end{equation}
Choosing $m=1$, $\hat K(\lambda )=X_{[-\pi ,\pi]}(\lambda )$ one obtains
the sine kernel $K(x-y)=\tfrac{\sin\pi (x-y)}{\pi (x-y)}$.  A special 
role played by the sine kernel can be highlighted by the fact that
$\tfrac{1}{\pi^2}\log L$ rate of the growth for Var$(\#_{[-L,L]})$ is the
slowest among all translation-invariant kernels $K(x-y)$ corresponding to 
projectors,
$\hat K = \chi_B $, for which $inf(B), \ \ sup(B)$ are the density points of 
$B$   (if $K$ is not a projector it follows from Lemma 6 that the rate of 
the growth of the variance is linear). 

As an example consider
$ B=\bigsqcup_{n\geq 1}[n,n+\tfrac{1}{n^\gamma}], \gamma>1$, than
one has $\tfrac{d\mu}{d\lambda}\sim\vert\lambda\vert^{1-\tfrac{1}{\gamma}}$ and
Var$(\#_{[-L,L]})\sim L^{\tfrac{1}{\gamma}}$.  More generally, $\tfrac{d\mu}
{d\lambda}\sim\vert\lambda\vert^\alpha ,0<\alpha <1$, implies Var$(
\#_{[-L,L]})\sim L^{1-\alpha}$.

\section{Central Limit Theorem for Counting Function and Empirical
Distribution Function of Spacings}

In [CL] Costin and Lebowitz proved the Central Limit Theorem for $\#_{[-L,L]}$
in the case of the sine kernel.  The article also contains a remark on 
p.~71, due to Widom, that the result holds for a larger class of Random
Matrix models.  In its general form this theorem appeared in [So2].
\medskip

\noindent{\bf Theorem 8}.  {\it Let $E$ be as in (1.1), $\{0<K_t\leq 1\}$ a
family of locally trace class operators in $L^2(E), \{(X,\b ,P_t)\}$
a family of the corresponding determinantal random point fields in $E$,
and $\{I_t\}$ a family of measurable subsets in $E$ such that
\begin{equation}
{\rm Var}_t\#_{I_t}=\tr (K_t\cdot \chi_{I_t}-(K_t\cdot \chi_{I_t})^2)\rightarrow
\infty\text{ \rm as }t\rightarrow\infty .
\end{equation}
Then the distribution of the normalized number of particles in $I_t$ (with
respect to $P_t$) converges to the normal law, i.e.,
$$\frac{\#_t-\bbE\#_{I_t}}{\sqrt{{\rm Var}_t\#_t}}
\overset{w}{\longrightarrow}N(0,1)$$}
\medskip

\noindent{\bf Remark 12}.  It was shown in [So2] that the condition (4.1)
from Theorem 8 (the growth of the variance) is satisfied for the Airy
kernel ($K_t\equiv K$ from (2.36), $I_t$ expanding), the Bessel kernel
($K_t\equiv K$ from (2.37), $I_t$ expanding) and for the families of
kernels $\{K_n\}$ corresponding to random matrices from the Classical
Compact Groups (\S 2.3b), \S 2.3c)).  In all these cases Var$_t\#_{I_t}$
growth logarithmically with respect to $\bbE_t\#_{I_t}$.
\medskip

\noindent{\bf Remark 13}.  To construct an example of the kernel $0\leq K\leq 
\id$ such that $E\#_{[-n,n]}=\tr K\cdot \chi_{[-n,n]}\rightarrow\infty$ as
$n\rightarrow\infty$, but Var$\#_{[-n,n]}=\tr (K\cdot \chi_{[-n,n]}-(K\cdot
\chi_{[-n,n]})^2)$ stays bounded, consider $\{\varphi_n(x)\}^\infty_{n=-\infty}$
satisfying

a) supp $\varphi_n\in (n,n+1)$,

b) $\Vert\varphi_n\Vert_{L^2}=1$.

Then $K(x,y)=\sum^\infty_{n=-\infty} (1-\tfrac{1}{n^2+1})\cdot\varphi_n(x)
\cdot\overline{\varphi_n(y)}$ is the desired kernel.  Indeed, $\bbE
\#_{[-n,n]}=\sum^n_{k=-n}(1-\tfrac{1}{k^2+1})\overset{n\rightarrow\infty}
{\longrightarrow}\infty ,{\rm Var} \#_{[-n,n]}=\sum^n_{k=-n}(1-\tfrac{1}
{k^2+1})\cdot\tfrac{1}{k^2+1}\rightarrow\sum^\infty_{-\infty}(1-\tfrac{1}
{k^2+1})\cdot\tfrac{1}{k^2+1}<\infty$.  From the other side if $0\leq
K\leq\id$ is compact, locally trace class  and $\tr K\cdot \chi_{[-n,n]}
\rightarrow +\infty$, then $\tr K\cdot \chi_{[-n,n]}-(K\cdot X_{[-n,n]})^2
\rightarrow +\infty$.

The result of Theorem 8 can be generalized to a finite number of intervals.
Namely, if $I^{(1)}_t,\dots ,I^{(m)}_t$ are disjoint subsets such that
Cov$_t (\#_{I^{(k)}_t}, \#_{I^{(j)}_t})/V_t\rightarrow b_{ij}$ as 
$t\rightarrow\infty , 1\leq i,j\leq m$, where $V_t$ is some function
of $t$ growing to infinity, then the distribution of $((\#_{I^{(k)}_t}
-\bbE_t\#_{I^{(k)}_t}/V_t^{\tfrac{1}{2}})_{k=1,\dots ,m}$ converges to
the $m$-dimensional centralized normal vector with the covariance matrix
$(b_{ij})_{1\leq i,j\leq m}$ (see [So2]).

Finally, we turn our attention to the problem of the global distribution
of spacings.  Let $E=\bbR^d$ or $\bbZ^d, \{B_j\}^k_{j=1}$ be some bounded
measurable subsets of $E$, and $\{n_j\}^k_{j=1}$ be some non-negative
integers.  We will be interested in the counting statistics of the
following type
\begin{equation}
\eta_L(B_1,\dots ,B_k;n_1,\dots ,n_k):=\#(x_i\in [-L,L]^d: \#_{x_i+B_j}
=n_j, j=1,\dots , k)
\end{equation}
We can assume without loss of generality that $\{B_j\}$ are disjoint and
do not include the origin.  If $d=1, k=1, B_1=(0,s]$, then $\eta_L
((0,s]),0)$ is the number of the nearest spacings in $[-L,L]$ greater than
$s: \eta_L((0,s],0)=\#\{x_i\in [-L,L]: x_{i+1}-x_i>s\}$, and $\eta_L
((0,s]),n)$ is the number of $n$-spacings greater than $s: \eta_L((0,s],n)=
\#\{x_i\in [-L,L]: x_{i+n+1}-x_i>s\}$.  In [So1] we proved the convergence
in law of the process $\tfrac{\eta_L((0,s],0)-\bbE\eta_L((0,s],0)}
{L^{\tfrac{1}{2}}}$ to the limiting Gaussian process in the case $K(x,y)=\tfrac
{\sin\pi (x-y)}{\pi (x-y)}$.  Recall that the convergence in law (Functional
Central Limit Theorem) implies not only the convergence of the 
finite-dimensional distributions, but also the convergence of functionals
continuous in the appropriate (e.g. locally uniform) topology on the space
of sample paths.  The proof of the Central Limit Theorem for the 
finite-dimensional distributions of $\eta_L((0,s],0)$ can be extended
essentially word by word to the case of arbitrary, not necessarily translation
invariant, kernel $K(x,y)$ and dimension $d\geq 1$, assuming the conditions
(4.33), (4.34), (4.35) are satisfied.  One can also replace $(0,s]$ by an
arbitrary measurable bounded $B\subset E$.  For the convenience of the
reader we sketch the main ideas of the proof of the finite-dimensional
Central Limit Theorem below.  Let us fix $B_1,\dots ,B_k; n_1,\dots ,n_k$.
We construct a new (called modified) random point field such that $\eta_L
(B_1,\dots B_k; n_1,\dots ,n_k)$ is equal to the number of all particles
of the modified random point field in $[-L,L]^d$.  Namely we keep only 
those particles of the original r.p.f. for which
\begin{equation}
\#_{x_i+B_j}=n_j, j=1,\dots ,k,
\end{equation}
and throw away the particles for which (4.3) is violated.  The modified
r.p.f. in general will no longer be a d.r.p.f.  What is important is that its
correlation functions and cluster functions (see Definition 6 below)
can be expressed in terms of the correlation functions of the original
determinantal r.p.f.  Let us denote by $\rho_\ell (x_1,\dots ,x_\ell ;
B_1,\dots ,B_k;n_1,\dots ,n_k)$ the $\ell$-point correlation function of
the modified r.p.f.  Suppose that
\begin{equation}
x_i\not\in x_j + B_p,\ 1\leq i\neq j\leq\ell ,\ 1\leq p\leq k.
\end{equation}
Then by the inclusion-exclusion principle
\begin{equation}
\begin{split}
&\rho_\ell (x_1,\dots ,x_\ell ;B_1,\dots ,B_k; n_1,\dots ,n_k)=
\sum^\infty_{m=0}\frac{(-1)^m}{m!}\\
&\underset{\longleftarrow\ell\text{ times} \longrightarrow}
{\int_{(x_1+B_1)^{n_1}\times\dots\times (x_1+B_k)^{n_k}}
\int_{(x_\ell +B_1)^{n_1}\times\dots\times (x_\ell +B_k)^{n_k}}}\\
&\int_{((x_1+\bigsqcup^k_{j=1}B_j)\bigsqcup\ldots\bigsqcup (x_\ell +
\bigsqcup^k_{j=1}B_j))^m}\rho_{\ell +\ell\cdot n+m}(x_1,\dots ,x_\ell ;\\
&x_{11},\dots ,x_{1n},x_{21},\dots ,x_{2n},\dots ,x_{\ell 1},\dots ,
x_{\ell n}, y_1,\dots y_m)dy_1\dots dy_m\\
& dx_{\ell 1}\dots dx_{\ell n}
\dots dx_{11}\dots dx_{1n},\\
& n=n_1+\dots +n_k.
\end{split}
\end{equation}
If (4.4) is violated then the formula is quite similar, the only difference
is that the exponent $n_j$ in $(x_i+B_j)^{n_j}=
\underset{\longleftarrow n_j\text{ times}\longrightarrow}
{(x_i+B_j)\times\dots\times(x_i+B_j)}$, $1\leq i\leq\ell , 1\leq j\leq k$,
has to be replaced by $n_j-\#(1\leq r\neq i\leq k: x_r\in x_i+B_j$).
While formulas (4.5) appear to be cumbersome and lengthy, they are
nevertheless quite useful for calculating the asymptotics of the moments
of $\eta_L(B_1,\dots ,B_k;n_1,\dots n_k)$. (Of course the assumption
that the correlation functions of the original r.p.f. are the determinants
is the key here.)  Recall the definition of the cluster functions.
\medskip

\noindent{\bf Definition 6}. {\it The $\ell$-point cluster functions 
$r_\ell (x_1,\dots ,x_\ell ), \ell =1,2,\dots ,$ of a random point field
are defined by the formula
\begin{equation}
r_\ell (x_1,\dots ,x_\ell )=\sum_G(-1)^{m-1}(m-1)!\cdot \prod^m_{j=1}
\rho_{\vert G_j\vert}(\bar x(G_j))
\end{equation}
where the sum is over all partitions $G$ of $[\ell ]=\{1,2,\dots ,\ell\}$
into subsets $G_1,\dots ,G_m,m=1,\dots ,\ell$, and $\bar x(G_j)=\{x_i:
i\in G_j\}, \vert G_j\vert =\#(G_j)$.}
\medskip

The cluster functions are also known in the Statistical Mechanics as the 
truncated correlated funciton and the Ursell functions.  Sometimes in 
the literature the r.h.s. of (4.6) defines $(-1)^{\ell -1}
r_\ell$.  Correlation functions can be obtained from cluster functions
by the inversion formula
\begin{equation}
\rho_\ell (x_1,\dots ,x_\ell )=\sum_G\prod^m_{j=1}r_{\vert G_j\vert}
(\bar x(G_j)).
\end{equation}
((4.6) is just the M\"obius inversion formula to (4.7).)
The integrals of\break $r_\ell (x_1,\dots ,x_\ell )$ over
$\underset{\longleftarrow\ell\text{ times}\longrightarrow}{[-L,L]^d\times
\dots\times [-L,L]^d}=[-L,L]^{\ell d}$ are closely related to the cumulants
$C_j(L)$ of the number of particles in $[-L,L]^d$ : 

$ V_1(L)=\int_{[-L,L]^d}
r_1(x_1)dx_1=C_1(L)=\bbE\#_{[-L,L]^d},$

$ V_2(L):=\int_{[-L,L]^d}
\int_{[-L,L]^d}r_2(x_1,x_2)dx_1\break dx_2=C_2(L)-C_1(L)={\rm Var}
\#_{[-L,L]^d}-\bbE\#_{[-L,L]^d}, $

  $V_3(L):=\int_{[-L,L]^d}\int_{[-L,L]^d}
\int_{[-L,L]^d}\break r_3(x_1,x_2,x_3)dx_1dx_2dx_3=C_3(L)-3C_2(L)+2C_1(L)$.

In general,
\begin{equation}
\sum^\infty_{n=1}\frac{C_n(L)}{n!}z^n=\sum^\infty_{n=1}\frac{V_n(L)}{n!}
(e^z-1)^n
\end{equation}
(see [CL], [So1]).  For the determinantal random point fields
\begin{equation}
r_\ell (x_1,\dots ,x_\ell )=(-1)^{\ell -1}\sum_{\text{cyclic }\sigma\in
S_\ell}K(x_1,x_2)\cdot K(x_2,x_3)\cdot\ldots\cdot K(x_\ell ,x_1),
\end{equation}
where the sum in (4.9) is over all cyclic permutations, and the term written
in the body of the sum corresponds to $\sigma =(1\ 2\ 3\dots\ell )$.  One
can also rewrite (4.9) as
\begin{equation}
\begin{split}
&r_\ell (x_1,\dots ,x_\ell )=(-1)^{\ell -1}\cdot\frac{1}{\ell}\sum_{\sigma
\in s_\ell}K(x_{\sigma (1)},x_{\sigma (2)})\cdot\\
&K(x_{\sigma (2)},x_{\sigma (3)})\ldots\cdot K(x_{\sigma (\ell )},x_{\sigma
(1)}).
\end{split}
\end{equation}
We note that the difference between (4.9) and the formula for $\ell$-point
correlation
\begin{equation}
\rho_\ell (x_1,\dots ,x_\ell )=\sum_{\sigma\in s_\ell}(-1)^\sigma K(x_1,
x_{\sigma (1)})\cdot K(x_2,x_{\sigma (2)})\cdot\ldots\cdot K(x,x_{\sigma
(\ell )})
\end{equation}
is that the summation in (4.9) is only over cyclic permutations.  It appears
that a relation between $\rho_\ell (x_1,\dots ,x_\ell ; B_1,\dots ,B_k;
n_1,\dots ,n_k)$ and $r_\ell (x_1,\dots ,x_\ell ;\break
B_1,\dots ,B_k; n_1,
\dots ,n_k)$ (at least when (3.19) is satisfied) is of a similar nature.
\medskip

\noindent{\bf Lemma 7}.  {\it Let (4.4) be satisfied.  Then
\begin{equation}
\begin{split}
&r_\ell (x_1,\dots ,x_\ell ;B_1,\dots ,B_k; n_1,\dots ,n_k)=
\sum^\infty_{m=0}\frac{(-1)^m}{m!}\\
&\underset{\longleftarrow\ell\text{ times}\longrightarrow}
{\int_{(x_1+B_1)^{n_1}\times\dots\times(x_1+B_k)^{n_k}}\dots
\int_{(x_\ell +B_1)^{n_1}\times\dots\times(x_\ell +B_k)^{n_k}}}\\
&\int_{((x_1+\bigsqcup^k_{j=1}B_j)\bigsqcup\dots\bigsqcup (x_\ell
+\bigsqcup^k_{j=1}B_j)^m}\rho_{\ell +\ell\cdot n+m,\ell}(x_1,\dots ,x_\ell ;
x_{11},\dots ,x_{1n},\\
& x_{21},\dots ,x_{2n}\dots ,x_{\ell 1},\dots ,x_{\ell n},y_1,\dots ,y_m)
dy_1\dots dy_m\\
&dx_{\ell 1}\dots dx_{\ell n}\dots dx_{11}\dots dx_{1n},
\end{split}
\end{equation}
where $\rho_{\ell +\ell\cdot n+m,\ell}$ is defined below in (4.13).}
\medskip

To define $\rho_{\ell +\ell\cdot n+m,\ell}$ recall that $\rho_{\ell +\ell
\cdot n+m}(x_1,\dots ,y_m)=\sum_{\sigma\in S_{\ell +\ell\cdot n+m}}(-1)^\sigma
\break K(x_1,
\sigma (x_1))\cdot\ldots\cdot K(y_m,\sigma (y_m))$, where $\sigma$
is a permutation on the set of variables $(x_1,\dots ,x_\ell , x_{11},
\dots ,x_{\ell n},y_1,\dots y_m)$.  We write
\begin{equation}
\begin{split}
&\rho_{\ell +\ell\cdot n+m,\ell}(x_1,\dots ,y_m)=\sideset{}{^*}\sum_{\sigma\in S_{\ell +
\ell\cdot n+m}}(-1)^\sigma K(x_1,\sigma (x_1))\cdot\ldots\cdot
K(y_m,\sigma (y_m)),
\end{split}
\end{equation}
where the summation in $\sum^*$ is over the permutations $\sigma$ satisfying
the following property:

Let $\tau$ be a multivalued map defined on $\{1,\dots ,\ell\}$ with the values in
$\{1,\dots ,\ell\}$:
\begin{equation}
\begin{split}
&\tau (i)=\{j:\sigma\biggl (\{x_i,x_{i1},\dots ,x_{in}\}\bigsqcup\biggl (\{
y_1,\dots ,y_m\}\bigcap\bigl (x_i+\bigsqcup^k_{p=1}B_p\bigr )\biggr )\biggr )\\
&\bigcap\biggl (\{x_j,x_{j1},\dots ,x_{jn}\}\bigsqcup\biggl (\{y_1,\dots,
y_m\}\bigcap\bigl (x_j+\bigsqcup^k_{p=1}B_p\bigr )\biggr )
\biggr )\neq\emptyset\};
\end{split}
\end{equation}
then for any $1\leq i,j\leq\ell$ there exists $N=N(i,j)$ such that
\begin{equation}
\tau^N(i)\ni j.
\end{equation}
\medskip

\noindent{\bf Remark 14}.  The proof of Lemma 7 in the case $d=1,
K(x,y)=\tfrac{\sin\pi (x-y)}{\pi (x-y)},\break B_1=(0,s], n_1=0$, was given
in \S 3 of [So1].  In the general case the argument is absolutely the same.
As a corollary of Lemma 7 we obtain
\medskip

\noindent{\bf Lemma 8}.  {\it Let 
\begin{equation}
\vert K(x,y)\vert\leq\psi (x-y),
\end{equation}
and (4.4) hold for the $\ell$-tuple $(x_1,\dots ,x_\ell )$.  Then for any
$\delta >0$ the following estimate takes place:
\begin{equation}
\begin{split}
&\vert r_\ell (x_1,\dots ,x_\ell ; B_1,\dots ,B_k; n_1,\dots ,n_k\vert\leq
\const (\ell ,\delta )\cdot\\
&\sum_{\text{cyclic }\sigma\in S_\ell}\biggl (\psi (x_2-x_1)\cdot\psi 
(x_3-x_2)\cdot\ldots\cdot\psi (x_1-x_\ell )\biggr )^{1-\delta}
\end{split}
\end{equation}}
\medskip

\noindent{$\Box$} For the proof of Lemma 8 we refer the reader to [So1] \S 3.
The key element of the proof is the upper bound on the absolute value of the
$m^{th}$ term in (4.12) by
\begin{equation*}
\begin{split}
&\const_1 (n,\ell )\cdot\frac{1}{m!}\const_2^m\cdot\min\biggl\{\const_3
(n,\ell ); (\ell +\ell n+m)!\cdot\\
&\sum_{\text{cyclic }\sigma\in S_\ell}\biggl (\psi (x_2-x_1)\cdot\psi 
(x_3-x_2)\cdot\ldots\cdot\psi (x_1-x_\ell)\biggr )\biggr\}.
\end{split}
\end{equation*}
\hfill\qed

\noindent
If $\psi^{1-\delta}\in L^2 (E)$ for some $0<\delta <1$, then $\int_{[-L,L]^d}
\dots \int_{[-L,L]^d}\psi (x_2-x_1)^{1-\delta}\cdot\ldots\cdot\psi 
(x_1-x_\ell )^{1-\delta}dx_1\dots dx_\ell\leq\const (\psi )\cdot\int_{[
-L,L]^d}\psi (x-y)^{2-2\delta}dxdy=O (L^d)$, therefore by Lemma 8
\begin{equation}
\begin{split}
&\int_{[-L,L]^{d\ell}\cap (4.4)}r_\ell (x_1,\dots ,x_\ell ;B_1,\dots ,
B_k; n_1,\dots n_k)\\
&dx_1\dots dx_\ell=O(L^d), \ell =1,2,\dots
\end{split}
\end{equation}
In particular
\begin{equation}
\begin{split}
&\bbE\eta_L(B_1,\dots ,B_k; n_1,\dots ,n_k)=V_1(L)=\\
& \int_{[-L,L]^d}r_1(x; B_1,\dots ,B_k; n_1,\dots ,n_k)dx=O
(L^d)
\end{split}
\end{equation}
Suppose that one could show
\begin{equation}
\begin{split}
&\var\ \eta_L(B_1,\dots ,B_k; n_1,\dots n_k)=V_1(L)+V_2(L)=\\
&\int_{[-L,L]^d}r_1(x; B_1,\dots B_k; n_1,\dots ,n_k)dx+\int_{[-L,L]^d}
\int_{[-L,L]^d}\\
&r_2(x_1,x_2;B_1,\dots ,B_k; n_1,\dots n_k)dx_1dx_2=\const\cdot L^d
(1+\bar o(1)),
\end{split}
\end{equation}

\begin{equation}
\begin{split}
&\int_{[-L,L]^{\ell d}\setminus (4.4)}r_\ell (x_1,\dots ,x_\ell ;
B_1,\dots ,B_k; n_1,\dots n_k)dx_1\dots dx_\ell\\
&=\bar o\left (L^{\frac{\ell d}{2}}\right ), \ell >2.
\end{split}
\end{equation}
Since the $\ell^{th}$ cumulant of $\eta_L$ is a linear combination of
$V_i(L), i=1,2,\dots \ell$ (see (4.8)), the estimates (4.18)--(4.21)
would imly that the $\ell^{th}$ cumulant of $\eta_L$ is $\const\cdot L
\cdot (1+\bar o(1))$ for $\ell =2$ and grows slower than 
$L^{\tfrac{\ell d}{2}}$ for $\ell >2$.  This in turn would imply that
while the second cumulant of $\tfrac{\eta_L-\bbE\eta_L}{\sqrt{\var\ \eta _L}}$
is 1, all the other cumulants of $\tfrac{\eta_L-\bbE\eta_L}
{\sqrt{\var\ \eta_L}}$ go to zero as $L\rightarrow +\infty$.  The last
statement is equivalent to the statement that the moments of $\tfrac{\eta_L-\bbE\eta_L}
{\sqrt{\var\ \eta_L}}$ converge to the moments of the normal distribution,
and in particular
$$\frac{\eta_L-\bbE\eta_L}{\sqrt{\var\ \eta_L}}\overset{w}
{\longrightarrow}N(0,1).$$
Of course the devil is in the details.  It turns out that there is no nice
extension of the formulas (4.12), (4.13) to the case when (4.4) is not
satisfied.  Below we show how one can overcome these difficulties in the case
of $\eta_L(B;0)$ (i.e. $k=1, n_1=0$).  We introduce the centralized 
$\ell$-point correlation functions by the formula
\begin{equation}
\rho^{(c)}_\ell(x_1,\dots ,x_\ell )=\sideset{}{^{**}}\sum_G\prod^m_{j=1}
r_{\vert G_j\vert}(\overline{x}(G_j),
\end{equation}
where $\sum^{**}$ is the sum over all partitions $G=\{G_1,\dots ,G_m\},
m=1,2,\dots$ of $\{1,\dots ,\ell\}$ into two- and more element subsets (i.e.
$\vert G_j\vert >1, j=1,\dots ,m)$.  It follows from (4.7), (4.22) that
\begin{equation}
\begin{split}
&\rho^{(c)}_\ell (x_1,\dots ,x_\ell )=\rho_\ell (x_1,\dots ,x_\ell )+
\sum^\ell_{p=1}(-1)^p\sum_{1\leq i_1<\dots <i_p\leq\ell}\prod^p_{s=1}\\
& \rho_1(x_{i_s})\cdot\rho_{\ell -p}\biggl ((x_1,\dots ,x_\ell)\setminus
(x_{i_1},\dots ,x_{i_p})\biggr )=\rho_\ell (x_1,\dots ,x_\ell )-\\
&\sum^\ell_{p=1}\sum_{1\leq i_1<\dots <i_p\leq\ell}\prod^p_{s=1}\rho_1
(x_{i_s})\cdot\rho^{(c)}_{\ell -p}\biggl ((x_1,\dots ,x_\ell )\setminus
(x_{i_1},\dots ,x_{i_p})\biggr ).
\end{split}
\end{equation}
Let us denote by $M^{(c)}_{(\ell )}(L)$ the integral of the centralized
$\ell$-point correlation function of the modified random point field over
$[-L,L]^{\ell d}$,
\begin{equation}
M^{(c)}_{(\ell )}(L)=\int_{[-L,L]^d}\dots \int_{[-L,L]^d}\rho^{(c)}_\ell
(x_1,\dots ,x_\ell ;B_1;0)dx_1\dots dx_\ell .
\end{equation}
We have
\begin{equation}
\begin{split}
&\sum^\infty_{\ell =0}\frac{t^\ell}{\ell !}\bbE(\eta_L-\bbE\eta_L)^\ell =
e^{-t\bbE\eta_L}\cdot\sum^\infty_{\ell =0}\frac{t^\ell}{\ell !}\bbE
\eta^\ell_L=e^{-t\bbE\eta_L}\cdot\sum^\infty_{\ell =0}\frac{(e^t-1)^\ell}
{\ell !}\\
&\bbE\eta_L\cdot
 (\eta_L-1)\cdot\ldots\cdot (\eta_L-\ell +1)=e^{-t\bbE\eta_L}
\cdot e^{(e^t-1)\bbE\eta_L}\cdot\\
&\sum^\infty_{\ell =0}\frac{(e^t-1)^\ell}
{\ell !}M^{(c)}_{(\ell )}(L)
\end{split}
\end{equation}
If we can show that
\begin{equation}
M^{(c)}_\ell (L)=\begin{cases}(2n-1)!!\cdot\const^n_1\cdot L^{nd}\cdot
(1+\bar o(1))\text{ for }\ell =2n,\\ 
\bar o\left (L^{\frac{\ell d}{2}}\right )\text{ for }\ell =2n+1,\end{cases}
\end{equation}
and
\begin{equation}
\bbE\eta_L=\const_2\cdot L^d\cdot (1+\bar o(1)),
\end{equation}
then (4.25) implies
\begin{equation}
\bbE (\eta_L-\bbE\eta_L)^\ell =\begin{cases}(2n-1)!!\cdot (\const_1+
\const_2)^n\cdot L^{nd}\cdot (1+\bar o(1))\text{ for }\ell =2n,\\
\bar o\left (L^{\frac{\ell d}{2}}\right )\text{ for }\ell =2n+1,
\end{cases}
\end{equation}
and
$$\frac{\eta_L-\bbE\eta_L}{L^{\frac{d}{2}}}\overset{w}{\longrightarrow}
N(0,\const_1+\const_2)$$
One can in principle calculate 
$M^{(c)}_\ell (L)$ from (4.12), (4.13).  Indeed, if
\begin{equation}
x_i-x_j\not\in B,
\end{equation}
(we remark that (4.29) is exactly (4.4) written in the case $k=1, n_1=0$),
then the expression for $\rho^{(c)}_\ell (x_1,\dots ,x_\ell ;B;0)$
can be obtained from (4.22), (4.12), (4.13).  Otherwise $\rho_\ell (x_1,
\dots ,x_\ell ;B; 0)=0$, and (4.23) implies
\begin{equation}
\begin{split}
&\rho^{(c)}_\ell (x_1,\dots ,x_\ell ;B;0)=\sum^\ell_{p=1}(-1)^p\cdot
\sum_{1\leq i_1<\dots <i_p\leq\ell}\prod^p_{s=1}\\
&r_1 (x_{i_s}; B;0)\cdot \rho_{\ell -p}\biggl ((x_1,\dots ,x_\ell )
\setminus (x_{i_1},\dots ,x_{i_p})\biggr ).
\end{split}
\end{equation}
If for an $(\ell -p)$-tuple $(x_1,\dots ,x_\ell )\setminus (x_{i_1},\dots
,x_{i_p})$ the condition (4.29) is not satisfied, then the corresponding
term $\rho_{\ell -p}((x_1,\dots ,x_\ell )\setminus (x_{i_1},\dots ,x_{i_p}))$
in (4.20) is zero.  If (4.29) is satisfied for $(x_1,\dots ,x_\ell )
\setminus (x_{i_1},\dots ,x_{i_p})$, then we iterate (4.23) again
\begin{equation*}
\begin{split}
& \rho_{\ell -p}\biggl ((x_1,\dots ,x_\ell )\setminus (x_{i_1},\dots ,
x_{i_p})\biggr )=\\
& \rho^{(c)}_{\ell -p}\biggl ((x_1,\dots ,x_\ell )\biggr )
\setminus (x_{i_1},\dots ,x_{i_p})\biggr )+\sum\dots
\end{split}
\end{equation*}
We claim
\medskip

\noindent{\bf Lemma 9}. {\it Let the condition (4.29) be not satisfied
for the $\ell$-tuple $(x_1,\dots ,x_\ell )$.  Then
\begin{equation}
\rho^{(c)}_\ell (x_1,\dots ,x_\ell ;B; 0)=\sum_{\emptyset\subseteq D\subset
\{1,\dots ,\ell\}} C_D\cdot\prod_{i\not\in D}r_1(x_i;B;0)\cdot
\rho^{(c)}_{\vert D\vert}(\bar x(D)),
\end{equation}
where
\begin{equation}
C_D=\sum_{\underset{\text{(4.29) is satisfied for }\bar x(A)}
{A\supseteq D,}}(-1)^{\vert A\vert}.
\end{equation}

\medskip

In particular $C_D=0$ if (4.29) is not satisfied for $\bar x(D)$  or if
there exists $1\leq i\leq\ell ,i\not\in D$, such that for any $1\leq j\leq
\ell\ x_i-x_j\not\in B\cup (-B)$.}

Proof easily follows from the above arguments.
\medskip

\noindent{\bf Theorem 9}.  {\it Let $(X,B,P)$ be a determinantal random
point field with the kernel
\begin{equation}
\vert K(x,y)\vert\leq\psi (x-y),
\end{equation}
where $\psi$ is a bounded non-negative function such that $\psi\cdot (\log
(\tfrac{\psi +1}{\psi}))^n\in L^2(E)$ for any $n>0$.  Let for $\eta_L(B;0)=\#
(x_i\in [-L,L]^d:\#(x_i+B)=0)$ we have 
\begin{equation}
\var\ \eta_L(B;0)=\sigma^2\cdot L^d\cdot (1+\bar 0(1))
\end{equation}
Then the Central Limit Theorem holds:
$$\frac{\eta_L(B,0)-\bbE\eta_L(B;0)}{L^{\frac{d}{2}}}\overset{w}
{\longrightarrow}N(0,\sigma^2).$$}
\medskip

\noindent{\bf Remark 15}.  If $\cov (\eta_L(B_i;0),\ \eta_L(B_j;0))=b_{ij}
\cdot L^d\cdot (1+\bar 0(1)),\ 1\leq i,j\leq p$, then 
\begin{equation}
\left (\frac{\eta_L(B_i;0)-\bbE\eta_L(B_i;0)}{L^{\frac{d}{2}}}\right 
)_{1\leq i\leq p}\overset{w}{\longrightarrow}N(0,(b_{ij})_{1\leq i,j\leq p}
).
\end{equation}
Recall that
\begin{equation*}
\begin{split}
&\cov (\eta_L(B_i;0); \eta_L(B_j;0))=\bbE (\eta_L(B_i; 0)-\bbE\eta_L (B_i,0))
\cdot (\eta_L(B_j;0)-\\
&\bbE\eta_L(B_j;0))=
\int_{\overset{[-L,L]^{2d}}{\cap\{x_1-x_2\not\in B_i\cup (-B_j)\}}}\bigl (
\sum^\infty_{m=0}\frac{(-1)^m}{m!}\cdot\int_{((x_1+B_i)\sqcup (x_2+B_j))^m}\\
&\rho_{2+m, 2}(x_1,x_2;y_1,\dots ,y_m)
dy_1\dots dy_m)dx_1dx_2-\int_{[-L,L]^d}r_1(x_1;B_i;0)\cdot\\
&\int_{(x_1+B_i)\cup (x_1-B_j)}r_1(x_2;B_j;0)dx_2dx_1+
\int_{[-L,L]^d}r_1(x;B_1\cup B_2,0)dx
\end{split}
\end{equation*}
\medskip

\noindent{\bf Remark 16}.  Lemma 8 suggests slightly more restrictive
condition on $\psi$, namely $\psi^{1-\delta}\in L^2(E)$ for some $0<\delta 
<1$.  However, looking at the proof of Lemma 6 one immediately realizes
that it is possible to replace $\psi^{1-\delta}$ in (4.17) by $\psi\cdot (
\log (\tfrac{\psi +1}{\psi}))^n$ with $n>3\ell$.
\medskip

\noindent{\bf Proof of Theorem 9}.  It follows from (4.25)-(4.28) that it
is enough to show
\begin{equation}
\begin{split}
&\int_{[-L,L]^{2nd}}\rho^{(c)}_{2n}(x_1,\dots ,x_{2n};B;0)dx_1\dots dx_{2n}=\\
&(2n-1)!!\cdot \biggl (\int_{\overset{[-L,L]^{2d}\cap}{\{x-y\not\in B
\cap (-B)\}}}r_2 (x,y;B;0)dxdy-\int_{[-L,L]^d}\\
&r_1(x;B;0)\int_{(x+B)\cup (x-B)}r_1(y;B;0)dydx\biggr )^n +\bar o(L^{nd}),\\ 
&n=1,2,\dots , 
\end{split}
\end{equation}

\begin{equation}
\begin{split}
&\int_{[-L,L]^{(2n+1)d}}\rho^{(c)}_{2n+1}(x_1,\dots ,x_{2n+1};B;0)dx_1\dots
dx_{2n+1}=\bar o\left (L^{\frac{2n+1}{2}d}\right ),\\
& n=1,2,\dots
\end{split}
\end{equation}
\medskip
\qed

\noindent{\bf Lemma 10}. {\it
\begin{equation}
\begin{split}
&\int_{[-L,L]^{2nd}\cap (4.29)}\rho^{(c)}_{2n}(x_1,\dots ,x_{2n};B;0)dx_1
\dots dx_{2n}=\\
&(2n-1)!!\cdot \biggl (\int_{\overset{[-L,L]^{2d}\cap}{\{x-y\not\in B\cup
(-B)\}}}r_2(x,y;B;0)dxdy\biggr )^n+\bar o(L^{nd}),
\end{split}
\end{equation}

\begin{equation}
\int_{[-L,L]^{(2n+1)d}\cap (4.29)}\rho^{(c)}_{2n+1}(x_1,\dots ,x_{2n+1};
B;0)dx_1\dots dx_{2n+1}=\bar o(L^{\frac{2n+1}{2}d}).
\end{equation}}

Recall that all $r_\ell (x_1,\dots ,x_\ell ;B;0)$ are bounded functions
(see (4.17)).  Let us rewrite (4.22) as
$$\rho^{(c)}_\ell (x_1,\dots ,x_\ell)=\sideset{}{'}\sum_G\prod^m_{j=1}
r_{\vert G_j\vert}(\bar x(G_j))+\sideset{}{''}\sum_G\prod^m_{j=1}
r_{\vert G_j\vert}(\bar x(G_j)),$$
where $\sum '$ is the sum over all partitions of $\{1,\dots ,\ell\}$ into
pairs, and $\sum ''$ is the sum over all other two- and more element
partitions.  Let $\ell$ be even, $\ell =2n$.  
Integrating $\sum_G'$ over $[-L,L]^{2nd}\cap (4.29)$ we obtain
exactly the r.h.s. of (4.38) (there are $(2n-1)!!$ partitions of $\{1,\dots
,2n\}$ into two-element sets).  It follows from (4.17) and the estimate
below Lemma 8 that $\int_{[-L,L]^{\ell d}}\vert r_\ell (x_1,\dots ,
x_\ell ;B;0)\vert dx_1\dots dx_\ell =\underline{O}(L^d)$.  Therefore the
integral of $\sum_G''$ over $[-L,L]^{2nd}\cap (4.29)$ is $\bar o(L^{nd})$.
The formula (4.39) can be proven in the same way.\hfill$\Box$

To estimate
\begin{equation}
\int_{[-L,L]^{2nd}\setminus (4.29)}\rho^{(c)}_{2n}(x_1,\dots ,x_{2n};
B;0)dx_1\dots dx_{2n}
\end{equation}
we introduce the equivalence relation on $\{x_1,\dots ,x_{2n}\}$ by calling
$x_i,x_j$ ``neighbors'' if there exists a sequence of indices $1\leq i_0,i_1,\dots 
,i_u\leq 2n,\ 1\leq u\leq 2n$, such that $i_0=i,\ i_u=j$,
and $x_{i_{s+1}}-x_{i_s}\in B\cup (-B),\ s=0,\dots ,u-1$.  We claim that
the contributions of order $O (L^{nd})$ appear in (4.40) only 
from such sets of $(x_1,\dots x_{2n})$ where each equivalence class of
``neighbors'' has either one or two indices.  Consider for example the case
when we have $k$ two-element classes $\{x_1,x_2\},\dots ,\{x_{2k-1},
x_{2k}\}$ and $2n-2k$ one-element equivalence classes $\{x_{2k+1}\},\dots
,\{x_{2n}\}$.  Similarly to the calculations on pp. 596-597 of [So1] we 
verify that the integral of $\rho^{(c)}_{2n}(x_1,\dots ,x_{2n};B;0)$ over
the subset of $[-L,L]^{2nd}$ corresponding to the above partition is equal
to
\begin{equation}
\begin{split}
&(2n-2k-1)!!\cdot\biggl (-\int_{[-L,L]^d}r_1(x;B;0)\int_{(x+B)\cup (x-B)}
r_1(y;B;0)dydx\biggr )^k\\
&\cdot\biggl (\int_{\underset{\{x-y\not\in B\cap (-B)\}}{[-L,L]^{2n}\cap}}
r_2(x,y;B;0)dxdy\biggr )^{n-k}+\bar o(L^{nd}).
\end{split}
\end{equation}
After the summation over all partitions into one- and two-element equivalence
classes of ``neighbors'' (we remark that (4.38) corresponds to the partition
into singletons), we obtain exactly (4.36).  It follows
from Lemma 7 and (4.17) that all other partitions into the equivalence
classes give negligible contributions. (4.37) can be proven in a similar
fashion.\hfill$\qed$
\medskip

The conditions of Theorem 9 are very unrestrictive in the case of translation
invariant kernels.  The covariance function of the limiting Gaussian process
w. $\lim\tfrac{\eta_L((0,\bar s],0)-\bbE\eta_L((0,\bar s];0)}
{L^{d/2}}$ is then given by the $d$-dimensional analogues of the
formulas (37), (38), (26) from [So1] (of course one has to replace $\tfrac
{\sin\pi (x-y)}{\pi (x-y)}$ by $K(x-y)$).  Here and below we denote by
$(0,\bar s]$ the rectangle $(0,s_1]\times\dots\times (0,s_d],\ \bar s=(
s_1,\dots ,s_d)$.  In particular, if $K(x)$ is continuously differentiable
the limiting Gaussian process is H\"older-continuous with any exponent less
than $\tfrac{1}{2}$.
Among other characteristics of the modified random point field (with respect
to $B=(0,\bar s],\ n=0)$ one may be interested in the spectral measure
of the restriction of the group $\{U^t\}$ to the
subspace of the centralized linear statistics.  We shall denote the spectral
measure by $\mu^{(s)}(d\lambda )$.  Recall that the spectral measure 
$\mu^{(0)}(d\lambda )=\mu (d\lambda )$ of the original determinantal random
point field is given by (3.12).  In particular, for the sine kernel
$$\frac{d\mu}{d\lambda}=\begin{cases}\frac{\vert\lambda\vert}{2\pi}, &
\vert\lambda\vert\leq 2\pi,\\1,&\vert\lambda\vert >2\pi.\end{cases}$$ 
After lengthy, but rather straightforward calculations one can  obtain that in
the case of the sine kernel :
\begin{equation}
\frac{d\mu^{(s)}}{d\lambda}=\frac{\pi^2s^3}{9}+\frac{\vert\lambda\vert}
{2\pi}\cdot\left (1-\frac{4}{3}\pi^2s^3\right )+O(s^4)+
O(\vert\lambda\vert\cdot s^4)+O (\vert\lambda^2
\vert\cdot s^2)
\end{equation}

We note that $\tfrac{d\mu^{(s)}}{d\lambda}(0)\neq 0$ if $s\neq 0$, $s$ small,
which is consistent with Var $\eta_L((0,s];0)\sim L$.  For the proof of
the Functional Central Limit Theorem we refer the reader to pp.~577,
598--600 of [So1].  Suppose that
\begin{equation}
L^{-d}\frac{\partial}{\partial s}\eta_L\bigl ((0,\bar s];0)\bigr ),L^{-d}\frac
{\partial}{\partial s}\cov (\eta_L\bigl ((0,\bar s];0\bigr );
\eta_L\bigl ((0,\bar t];0\bigr )\bigr )
\end{equation}
are uniformly bounded in $L,\bar s,\bar t$, where $\bar s,\bar t$ belong to
compact subsets of $\bbR^d_+(\bbZ^d_+)$.  By smoothing with a $C^\infty$
approximate $\delta$-function one can construct a continuous approximation
$\tilde\eta_L((0,\bar s];0)$ such that $\vert\tilde\eta_L((0,\bar s];0)
-\eta_L((0,\bar s];0)\vert\leq 1$.  As a result
$$\frac{\tilde\eta_L((0,\bar s];0)-\bbE\tilde\eta_L((0,\bar s];0)}
{L^{\frac{d}{2}}}$$
is a random continuous function in $\bar s$, and
\begin{equation}
\biggl\vert\frac{\tilde\eta_L((0,\bar s];0)-\bbE\tilde\eta_L((0,\bar s];0)}
{L^{\frac{d}{2}}}-\frac{\eta_L((0,\bar s];0)-\bbE\eta_L((0,\bar s];0)}
{L^{\frac{d}{2}}}\biggr\vert\leq\frac{2}{L^{\frac{d}{2}}}
\end{equation}
The distribution of the random process
$\frac{\tilde\eta_L((0,\bar s];0)-\bbE\tilde\eta_L((0,\bar s];0)}
{L^{\frac{d}{2}}}$ defines a probability measure on $C([0,\infty )^d)$. 
By the convergence in law of random processes we mean the weak convergence
of the induced probability measures on $C([0,\infty )^d)$ (see [B],
in general one can consider different spaces of sample paths, e.g. the
space of c\`adl\`ag functions, instead of the space of continuous functions).
\medskip

\noindent{\bf Theorem 10}.  {\it Let the condition (4.33), (4.34), (4.35),
(4.43) be satisfied. Then the random process
$$\frac{\tilde\eta_L((0,\bar s];0)-\bbE\tilde\eta_L((0,\bar S];0)}
{L^{\frac{1}{2}}}$$
converges in law to the limiting Gaussian process.}

\def\am{{\it Ann. of Math.} }
\def\ap{{\it Ann. Probab.} }
\def\temf{{\it Teor. Mat. Fiz.} }
\def\jmp{{\it J. Math. Phys.} }
\def\cmp{{\it Commun. Math. Phys.} }
\def\jsp{{\it J. Stat. Phys.} }
\def\prl{{\it Phys. Rev. Lett.} }
\def\arma{{\it Arch. Rational Mech. Anal.} }

\end{document}